\documentstyle[12pt]{article} 
\begin{document}
\newcommand{\singlespace}{
    \renewcommand{\baselinestretch}{1}
\large\normalsize}
\newcommand{\doublespace}{
   \renewcommand{\baselinestretch}{1.2}
   \large\normalsize}
\renewcommand{\theequation}{\thesection.\arabic{equation}}

\input amssym.def
\input amssym
\setcounter{equation}{0}
\def \ten#1{_{{}_{\scriptstyle#1}}}
\def \Z{\Bbb Z}
\def \C{\Bbb C}
\def \R{\Bbb R}
\def \Q{\Bbb Q}
\def \N{\Bbb N}
\def \l{\lambda}
\def \V{V^{\natural}}
\def \wt{{\rm wt}}
\def \tr{{\rm tr}}
\def \Res{{\rm Res}}
\def \End{{\rm End}}
\def \Aut{{\rm Aut}}
\def \mod{{\rm mod}}
\def \Hom{{\rm Hom}}
\def \im{{\rm im}}
\def \<{\langle} 
\def \>{\rangle} 
\def \w{\omega}
\def \o{\omega}
\def \t{\tau }
\def \char{{\rm char}}
\def \a{\alpha }
\def \b{\beta}
\def \e{\epsilon }
\def \la{\lambda }
\def \om{\omega }
\def \O{\Omega}
\def \qed{\mbox{ $\square$}}
\def \pf{\noindent {\bf Proof: \,}}
\def \voa{vertex operator algebra\ }
\def \voas{vertex operator algebras\ }
\def \p{\partial}
\def \1{{\bf 1}}
\singlespace
\newtheorem{thmm}{Theorem}
\newtheorem{th}{Theorem}[section]
\newtheorem{prop}[th]{Proposition}
\newtheorem{coro}[th]{Corollary}
\newtheorem{lem}[th]{Lemma}
\newtheorem{rem}[th]{Remark}
\newtheorem{de}[th]{Definition}
\newtheorem{slem}[th]{Sublemma}

\begin{center}
{\Large {\bf Modularity in orbifold theory for vertex operator superalgebras}} \\
\vspace{0.5cm}

Chongying Dong\footnote{Supported by NSF grants, a China NSF grant and  a Faculty
research grant from  the University of 
California at 
Santa Cruz.} and Zhongping Zhao
\\
Department of Mathematics, University of
California, Santa Cruz, CA 95064 \\
\end{center}
\hspace{1.5 cm}
\begin{abstract} This paper is about the orbifold theory for vertex operator superalgebras.
Given a vertex operator superalgebra $V$ and a finite automorphism group $G$ of $V,$
we show that the trace functions associated to the twisted sectors are holomorphic in the
upper half plane for any commuting pairs in $G$ under the $C_2$-cofinite condition. We also 
establish that these functions afford a representation of the full modular group  if $V$
is $C_2$-cofinite and $g$-rational for any $g\in G.$  
\end{abstract}

\section{Introduction}

Modular invariance of trace functions is fundamental in the theory
of vertex operator algebras and conformal field theory (see for example
\cite{B2}, \cite{DVVV}, \cite{DLM3}, \cite{FLM3},\cite{H},\cite{Hu}, \cite{M2}, \cite{MS}, \cite{V}, \cite{Z}). Many 
important results such as the existence of twisted sectors in orbifold
theory, the Lie algebra structure of weight one subspace of a good vertex 
operator algebra, and classification of holomorphic vertex operator
algebras with small central charges can be obtained by studying the modular
invariance property of 1-point functions on torus (see \cite{DM1}, \cite{DM2},
\cite{Mo}, \cite{S}). The modular invariance of trace functions 
associated to modules for an arbitrary rational vertex operator algebra was 
first established under the $C_2$-cofinite condition in \cite{Z}. This
work was extended in \cite{DLM3} to include a finite automorphism
group of the vertex operator algebra. There were further generalizations in 
\cite{M1}, \cite{M2}, \cite{Y}.

Partially motivated by the generalized moonshine conjecture \cite{N}, the
trace functions associated to twisted modules for a vertex operator
algebra $V$ and a finite automorphism group $G$ was studied in \cite{DLM3}.
Under some finiteness conditions on a  rational vertex operator algebra 
$V,$   among other things, the precise number of inequivalent, simple 
$g$-twisted $V$-modules was determined and the modular invariance
(in a suitable sense) of the 1-point functions on torus associated
to any commuting pairs in $G$ was established. If $V$ is holomorphic,
the existence and the uniqueness of the twisted sector for each 
$g\in G$ was obtained. The important case in which $V=V^{\natural}$
is the   moonshine vertex operator algebra \cite{FLM3} and $G$ is 
the monster simple group \cite{G} plays a special
role in a proof  of 
several parts of the generalized moonshine conjecture \cite{DLM3} except 
the genus zero property. Also see \cite{T1},\cite{T2} and \cite{T3}
on the generalized moonshine conjecture.

In this paper we present a general theory of modular invariance
of trace functions associated to twisted sectors for a vertex
operator superalgebra $V$ and a finite automorphism group $G.$
Note that $V$ has a  canonical central automorphism 
$\sigma$ arising from the superspace structure of $V,$ which is crucial in formulating the main results in this paper.
Let $\bar G$ be the group generated by $G$ and $\sigma.$  
We now explain the main results in details.

Let $g\in G$ have order $T$ and $\sigma g $ have order $T'.$ Then
a simple $\sigma g$-twisted module  $M=(M,Y_M)$ has a grading of the shape
\begin{equation}\label{e1.3} 
M=\bigoplus_{n=0}^{\infty}M_{\lambda+ n/T}
\end{equation}
with $M_{\l}\ne 0$ 
for some complex number $\l$ which is called the {\em conformal weight}
 of $M$ (see Section 4 and \cite{DZ}). One basic result in \cite{DZ} is that
if $V$ is $g$-rational there are only finitely many simple $\sigma g$-twisted
$V$-modules up to isomorphism. As in \cite{DLM3},  for any $h\in \bar G$
we can define a $h\sigma gh^{-1}$-twisted module $h\circ M$ so that
$h\circ M=M$ as vector spaces and $Y_{h\circ M}(v,z)=Y_M(h^{-1}v,z)$ 
(see Section 6). $M$ is called $h$-stable if $M$ and $h\circ M$ are isomorphic.
As a result the stabilizer of $M$ in $\bar G$  acts projectively on
$M,$ denoting by $\phi.$  

Let $v\in V$ and  $M$ be $\sigma h$-stable. We set 
\begin{equation}\label{e1.8}
T_M(v, (g, h),q)  =  q^{\lambda-c/24}\sum_{n=0}^{\infty}\tr_{M_{\lambda+n/T}}
(o(v) \phi(\sigma h))q^{n/T}
\end{equation}
where $Y_M(v,z)=\sum_{m\in \frac{1}{T'}\Z}v(m)z^{-m-1}$ and 
$o(v)=v(\wt v-1)$  induces a linear map on each homogeneous subspace of $M.$
We extend the notation $T_M(v,(g,h),q)$ to all $v\in V$ linearly. 
 
As in \cite{Z}, there is another vertex operator superalgebra structure
on $V$ (see Section 3). If the original vertex operator superalgebra
is defined over sphere, the second vertex operator superalgebra 
can be regarded as defined over a torus. 
 We will denote the corresponding weight by $\wt [v].$
Here are our main results on the modular invariance.

\begin{thmm}\label{thmm3} Suppose that $V$ is a $C_2$-cofinite 
vertex operator superalgebra and $G$ a finite group of automorphisms of $V.$

  (i) Let $g,h\in G$ and $M$ be a simple, $\sigma h$-stable,
 $\sigma g$-twisted $V$-module.  Then 
the trace function  $T_M(v,(g,h),q)$ converges
to a holomorphic function in the upper half plane ${\frak h}$ 
where $q=e^{2\pi i\tau}$ and $\tau\in \frak h.$

  (ii) Suppose in addition that $V$ is $x$-rational for each $x\in \bar G.$ 
Let $v\in V$ satisfy $\wt[v] = k.$ Then the space of (holomorphic) functions 
in ${\frak h}$ spanned
by the trace functions $T_M(v, (g,h), \tau)$ for all choices of $g, h$ in $G$
 and 
$\sigma h$-stable $M$
 is  a (finite-dimensional)  $SL(2,\Z)$-module such that
 $$
T_M|\gamma(v, (g, h), \tau)  =  (c\tau+ d)^{-k} T_M (v, (g, h),  \gamma\tau),$$
where $\gamma\in SL(2,\Z)$ acts on $\frak h$ as usual.

 More precisely, if $\gamma=\left(\begin{array}{cc}
a & b\\ c& d\end{array}\right)\in SL(2,\Z)$ then we have an equality
$$
T_M(v, (g, h),\frac{a\tau+b}{c\tau+d})=(c\tau+d)^k\sum_{W}\gamma_{M,W}T_W(v,(g^ah^c,g^bh^d),\tau),$$
where $W$ ranges over the $\sigma g^a h^c$-twisted sectors which are $g^b h^d$ 
and $\sigma$-stable. The constants $\gamma_{M,W}$ depend only on $M, W$ and 
$\gamma$ only.
\end{thmm}

\begin{thmm}\label{thmm4} Let $V$ be  a $C_2$-cofinite and $\sigma$-rational
vertex operator superalgebra. Let $\{M^1,...,M^n\}$ be the inequivalent
$\sigma$-twisted $V$-modules which are $\sigma$-stable. Then for any
$\gamma=\left(\begin{array}{cc}
a & b\\ c& d\end{array}\right)$ in  $SL(2,\Z)$ there 
are constants $\gamma_{ij}$ for $1\leq i,j\leq n$ such that 
for any $v\in V$ with $\wt[v] = k,$ 
$$T_{M^i}(v, (1, 1),\frac{a\tau+b}{c\tau+d})=(c\tau+d)^k\sum_{j=1}^n\gamma_{ij}T_{M^j}(v,(1,1),\tau).$$
\end{thmm}

It is important to point out that Theorem \ref{thmm4} is an analogue 
of the main theorem in \cite{Z}. In fact, if $V_{\bar 1}=0$ then $V$
is a vertex operator algebra and $\sigma=id_V.$ In this  case, Theorem 
\ref{thmm4} reduces the main theorem in \cite{Z}. In particular,
the space of characters of irreducible $V$-modules is invariant 
under the action of the modular group. But
if $V_{\bar 1}\ne 0$ the space of characters of irreducible $V$-modules
is no longer invariant under the action of the modular group. However,
the space of $\sigma$-trace on $\sigma$-twisted sectors is modular
invariant. So the automorphism $\sigma$ plays a fundamental role in the study 
of modular invariance for vertex operator superalgebra and
its orbifold theory. This will be illustrated in Section 9 with examples.

There is an assumption in both Theorems \ref{thmm3} and \ref{thmm4}
on twisted modules. That is, the twisted modules are required 
to be $\sigma$-stable. There are examples in which irreducible $\sigma$-twisted
modules are not $\sigma$-stable. So these modules which are
not $\sigma$-stable are excluded in the discussion of modular invariance. 
To include these modules, one can use the Theta group $\Gamma_{\theta}
=\Gamma(2)\cup \Gamma(2)S$ instead of
the full modular group $SL(2,\Z).$ Here $\Gamma(2)$ is the congruence 
subgroup of $SL(2,\Z)$ consisting of elements which are congruent 
to the identity modulo 2 and 
$S=\left(\begin{array}{cc}0 &-1\\ 1 &0\end{array}\right).$  In fact if $G=1$ this has been done in \cite{H}. The $\sigma$ does not play any role in this
approach.

Again Let $G$ be a finite automorphism group of $V$ and $g\in G.$ Let
$M$ be a simple $g$-twisted $V$-module.  Let $G_M$ be the stablizer of
$M$ in $G.$ Then $M$ is $h$-stable for any $h\in G_M$ and $G_M$ 
acts projectively on $M.$ As before we use $\phi$ to denote the
projective representation. For $h\in G_M$ and $v\in V$ set
$$F_M(v, (g, h),q)  = \tr_M o(v)\phi(h)q^{L(0)-c/24}.$$

The analogue of Theorem \ref{thmm3} is the following which extends 
the modular invariance result in \cite{H} to the orbifold theory.
\begin{thmm}\label{thmm5} Suppose that $V$ is a $C_2$-cofinite 
vertex operator superalgebra and $G$ a finite group of automorphisms of $V.$

  (i) Let $g,h\in G$ and $M$ be a simple, $h$-stable  $g$-twisted $V$-module.  Then 
the trace function  $F_M(v,(g,h),q)$ converges
to a holomorphic function in the upper half plane ${\frak h}.$ 

  (ii) Suppose in addition that $V$ is $x$-rational for each $x\in G.$ 
Let $v\in V$ satisfy $\wt[v] = k.$ Then the space of (holomorphic) functions 
in ${\frak h}$ spanned
by the trace functions $F_M(v, (g,h), \tau)$ for all choices of $g, h$ in $G$
 and 
$h$-stable $M$  is  a (finite-dimensional)  $\Gamma_{\theta}$-module. That is,
if $\gamma=\left(\begin{array}{cc}
a & b\\ c& d\end{array}\right)\in \Gamma_{\theta}$ then we have an equality
$$F_M(v, (g, h),\frac{a\tau+b}{c\tau+d})=(c\tau+d)^k\sum_{W}\gamma_{M,W}F_W(v,(g^ah^c,g^bh^d),\tau),$$
where $W$ ranges over the $g^a h^c$-twisted sectors which are $g^b h^d$-stable. 
The constants $\gamma_{M,W}$ depend only on $M, W$ and 
$\gamma$ only.
\end{thmm}

This paper is a super version of the theories developed in \cite{DLM3}. The
main ideas and the broad outline of our proof follow from those in \cite{DLM3},
\cite{Z} and \cite{H}. Two important results on $g$-rationality and
the associative algebras $A_g(V)$ in \cite{DZ} (also see Section 4)
also play basic roles in this paper. If $V$ is $g$-rational, then
$A_g(V)$ is semisimple and 
the category of $V$-modules and the category of finite dimensional
$A_g(V)$-modules are equivalent. Furthermore there
are finitely many irreducible $g$-twisted modules. These results enable us
to prove that the space of 1-point $( g, h)$ functions is finite dimensional.
In fact, from the proof of the main results in this paper one can
see that we only need the assumptions that $V$ is $C_2$-cofinite and
$A_g(V)$ is a semisimple associative algebra. Although there are a lot
of expectations on the relation among rationality, $C_2$-cofiniteness,
the semisimplicity of $A_g(V),$ we cannot remove either 
the $C_2$-cofinite condition or the $A_g(V)$ semisimplicity condition. 

This paper is organized as follows. In Section 2 we review several family
of modular functions and elliptic functions, and their transformation
laws following \cite{DLM3}. Section 3 is about the definition
of vertex operator superalgebra and vertex operator superalgebra
on torus. In Section 3 we define weak, admissible and ordinary
$g$-twisted module for a vertex operator superalgebra $V$ 
and a finite order automorphism $g.$ We also present 
the main results concerning the associative algebra $A_g(V)$
and $g$-rationality.  We define the space of abstract 
1-point $(g,h)$-functions in Section 5 and establish that
the elements in the full modular group transform the space
of 1-point $(g,h)$ functions to  another space of 1-point functions.
We also discuss the shapes of $q$-expansions of the 1-point functions and
prove that these functions satisfy differential equations with regular 
singularity. This leads us to the definition of formal 1-point $(g,h)$
functions. Sections 6 and 8 are the key sections in this paper.
We show in Section 6 that the $\sigma h$-trace on $\sigma g$-twisted
modules produce 1-point $(g,h)$-functions. In fact, the inequivalent
simple $\sigma g$-twisted modules give linearly independent 
1 point functions. In Section 8 we show that if $V$ is $\sigma g$-rational
then the space of 1-point $(g,h)$ functions has a basis 
consisting of the $\sigma h$-trace on the inequivalent simple 
$\sigma g$-twisted modules. We also give several important corollaries 
in Section 8 concerning the rationality of the central charges and
the conformal weights of simple $g$-twisted modules. In Section 9 we
present the modular invariance result on $\Gamma_{\theta}.$ 
We give an existence result on the twisted sectors in Section 7 and discuss 
several
examples in Section 10 to illustrate the main theorems.

\section{$P$-functions and $Q$-functions}
\setcounter{equation}{0}
In this section , we review some properties of  the $P-$functions and $Q-$ functions
following \cite{DLM3}. These functions will be used extensively in later sections. 

Recall that the {\em modular group} {\bf $\Gamma=SL(2,{\Bbb Z})$} acts on the upper half plane
$\frak h=\{z\in {\Bbb C}|im z>0 \}$ via M\"obius 
transformations
\begin{equation}\label{g4.1}
\left(\begin{array}{cc}
a & b\\
c & d
\end{array}
\right)\tau=\frac{a\tau+b}{c\tau+d}.
\end{equation} 
$\Gamma$ also acts on the right of $S^1\times S^1$ via
\begin{equation}\label{g4.2}
(\mu,\lambda)
\left(
\begin{array}{cc}
a & b\\
c & d
\end{array}
\right)
=(\mu^a\lambda^c,\mu^b\lambda^d).
\end{equation}

Let   $\mu=e^{2\pi ij/M}$ and $\lambda=e^{2\pi il/N}$ for 
integers $j,l,M,N$ with $M,N>0.$ For each integer $k=1,2,\cdots$ and
each $(\mu,\lambda)$ we define a function $P_k$ on $\C\times \frak h$ as follows
 \begin{equation}\label{g1.3}
P_k(\mu,\lambda,z,q_{\tau})=P_k(\mu,\lambda,z,\tau)=
\frac{1}{(k-1)!}\sum^{\ \ \ \ \ \prime}_{n\in \frac{j}{M}+\Z}\frac{n^{k-1}q_z^n}{1-\l q_{\tau}^n}
\end{equation}
 In the above equation, $q_x=e^{2\pi ix}$ and  $\sum'$ means omit the term $n=0$ when $(\mu,\l)=(1,1).$ Then $P_k$ converges uniformly and absolutely on compact subsets of the 
region $|q_{\tau}|<|q_z|<1.$ We have the following transformation properties (see
Theorem 4.2 of \cite{DLM3}). 

\begin{th}\label{t1.1} Suppose that $(\mu,\l)\ne (1,1).$ Then  the following equalities hold:
$$P_k(\mu,\lambda,\frac{z}{c\tau+d},\gamma\tau)=
(c\tau+d)^kP_k((\mu,\lambda)\gamma,z,\tau).$$
for all $\gamma=
\left(
\begin{array}{cc}
a & b\\
c & d
\end{array}
\right)\in \Gamma.$
\end{th}

In order to define $Q-$functions we need to recall the Bernoulli polynomials $B_r(x)\in \Q[x]$
which are determined by  
$$\frac{te^{tx}}{(e^t-1)}=\sum^{\infty}_{r=0}\frac{B_r(x)t^r}{r!}.$$
For example, $B_0(x)=1,B_1(x)=x-\frac{1}{2},B_2(x)=x^2-x+\frac{1}{6}.$

For $(u,\lambda)=(e^{\frac{2 \pi ij}{M}},e^{\frac{2 \pi il}{N}})$ and $(u,\lambda)\neq (1,1),$when $k\geq 1$ and $k\in {\Bbb Z},$ we define 
\begin{eqnarray}
& &Q_{k}(\mu,\la,q_{\t})=Q_k(\mu,\l,\tau)\nonumber\\
& &\ \ \ \ \ \ =\frac{1}{(k-1)!}\sum_{n\geq 0}
\frac{\la(n+j/M)^{k-1}q_{\tau}^{n+j/M}}{1-\lambda q_{\t}^{n+j/M}}\nonumber\\
& &\ \ \ \ \ \ \ \ +\frac{(-1)^k}{(k-1)!}\sum_{n\geq 1}\frac{\la^{-1}(n-j/M)^{k-1}
q^{n-j/M}_{\tau}}{1-\lambda^{-1}q^{n-j/M}_{\tau}}-\frac{B_k(j/M)}{k!}\label{m4.23}.
\end{eqnarray}
Here $(n+j/M)^{k-1}=1$ if $n=0, j=0$ and $k=1.$ Similarly,
$(n-j/M)^{k-1}=1$ if $n=1, j=M$ and $k=1.$  
We also define 
\begin{equation}\label{g1.5}
Q_0(\mu,\l,\tau)=-1.
\end{equation}

Here are the transformation properties of $Q$-functions (see Theorem 4.6
of \cite{DLM3}).
\begin{th}\label{t1.2} If $k\geq 0$ then $Q_k(\mu,\l,\t)$  is a holomorphic 
modular form of weight $k.$ If $\gamma=\left(\begin{array}{cc} a & b\\ c & d
\end{array}\right)\in \Gamma$ it satisfies 
$$Q_k(\mu,\la,\gamma\tau)=(c\tau+d)^kQ_k((\mu,\la)\gamma,\tau).$$
\end{th}

For the further discussion we also define $\bar P-$functions.
Again we take $\mu=e^{2\pi ij/M}$ and
$\l=e^{2\pi il/N}.$ Set for $k\geq 1,$
\begin{equation}\label{m4.29}
\bar P_{k}(\mu,\la,z,q_{\t})=\bar P_{k}(\mu,\la,z,\tau)=\frac{1}{(k-1)!}\sum_{n\in j/M+\Z}^{\ \ \ \ \ \prime}\frac{n^{k-1}z^n}{1-\la q^n_{\tau}}
\end{equation}
with
\begin{equation}\label{m4.30}
\bar P_0=0.
\end{equation} Then we have (see Proposition 4.9 of \cite{DLM3}):
\begin{prop}\label{p4.8}
If $m\in\Z,$ $\mu=e^{2\pi ij/M}, k\geq 0$ and $(\mu,\la)\ne(1,1),$ then
\begin{eqnarray*} 
& &\ \ Q_{k}(\mu,\la,\t)+\frac{1}{k!}B_k(1-m+j/M)\\
& &=\Res_{z}\left(\sum_{s\geq 0}\frac{z_1^s}{z^{s+1}}z_1^{m-j/M}z^{-m+j/M}\bar P_{k}(\mu,\la,\frac{z_1}{z},\t)\right)\nonumber\\
& &\ \ +\Res_{z}\left(\la\sum_{s\geq 0}\frac{z^s}{z_1^{s+1}}z_1^{m-j/M}z^{-m+j/M}\bar P_{k}(\mu,\la,\frac{z_1q_{\tau}}{z},\tau)\right).
\end{eqnarray*}
\end{prop}

We also recall the Eisenstein series of weight $k$ for  even $k\geq 3:$
\begin{eqnarray}\label{m4.25}
& &G_k(\tau)=\sum^{\ \ \ \ \ \prime}_{m_1,m_2\in\Z}\frac{1}{(m_1\tau+m_2)^k}\nonumber\\
& &\ \ \ \ =(2\pi i)^k\left(\frac{-B_k(0)}{k!}+\frac{2}{(k-1)!}\sum_{n=1}^{\infty}\sigma_{k-1}(n)q^n\right)
\end{eqnarray}
where $\sigma_{k-1}(n)=\sum_{d|n}d^{k-1}.$ For this range of values
of $k$, $G_k(\t)$ is a modular form on $SL(2,\Z)$ of weight $k.$
We will  use of the normalized Eisenstein series 
\begin{equation}\label{m4.26}
E_k(\tau)=\frac{1}{(2\pi i)^k}G_k(\tau)=\frac{-B_k(0)}{k!}+\frac{2}{(k-1)!}\sum_{n=1}^{\infty}\sigma_{k-1}(n)q^n
\end{equation}
for even $k\geq 2.$ Note that $E_2(\tau)$ is not a modular form but it is a well defined 
holomorphic function in $0<|q|<1.$

\section{Vertex operator superalgebras}
\setcounter{equation}{0}

In this section we recall the definition of vertex operator superalgebra
(cf. \cite{B1}, \cite{DL}, \cite{FLM3}). We also study another vertex
operator superalgebra structure on torus following \cite{Z} to 
give the modularity result in later sections.

Let $V=V_{\bar 0}\oplus V_{\bar 1}$ be any ${\Bbb Z}_2$-graded vector space.
The element in $V_{\bar 0} $ (resp. $V_{\bar 1}$) is called  even (resp. odd). 
For any $v\in V_{\bar i}$ with  $i=0,1$ we define ${\tilde v}=i.$
For convenience, we let $\epsilon_{u,v}=(-1)^{{\tilde u}{\tilde v}}$ and $\epsilon_{v}=(-1)^{\tilde v}.$

\begin{de}{\rm A vertex operator superalgebra }(VOSA) is a quadruple $(V,Y,1,\o),$ where $V=V_{\bar 0}\oplus V_{\bar 1}$ is a ${\Bbb Z}_2$-graded vector space, ${\bf 1}\in V_0$ is the {\rm vacuum }
of V, $\o\in V_2$ is called {\rm conformal vector} of $V,$   and $Y$ is a linear map 
\begin{eqnarray*}\label{0a3}
& & V \to (\End\,V)[[z,z^{-1}]] ,\\
& & v\mapsto Y(v,z)=\sum_{n\in{\Z}}v(n)z^{-n-1}\ \ \ \  (v(n)\in
\End\,V).\nonumber
\end{eqnarray*}
satisfying the following axioms for $u,v\in V:$ \\
(i) For any $u,v\in V,$ $u(n)v=0$ for sufficiently large $n,$\\
(ii) $Y({\bf 1},z)=Id_{V},$\\
(iii) $Y(v,z){\bf 1}=v+\sum_{n\geq 2}v(-n){\bf 1}z^{n-1},$\\
(iv) The component operators of  $Y(\w,z)=\sum_{n\in\Z}L(n)z^{-n-2}$ satisfy the Virasoro algebra 
relation with {\em central charge} $c\in \C:$ 
\begin{equation}\label{g2.1}
[L(m),L(n)]=(m-n)L(m+n)+\frac{1}{12}(m^3-m)\delta_{m+n,0}c,
\end{equation}
and 
\begin{eqnarray}
& &L(0)|_{V_n}=n,\ n\in\Z \label{g2.1'}\\
& &\frac{d}{dz}Y(v,z)=Y(L(-1)v,z),\label{g2.3}
\end{eqnarray}
(v) The {\em Jacobi identity} for ${\Bbb Z}_2$-homogeneous $u,v\in V$ holds,  
\begin{equation}\label{2.4}
\begin{array}{c}
\displaystyle{z^{-1}_0\delta\left(\frac{z_1-z_2}{z_0}\right)
Y(u,z_1)Y(v,z_2)-{\epsilon}_{u,v} z^{-1}_0\delta\left(\frac{z_2-z_1}{-z_0}\right)
Y(v,z_2)Y(u,z_1)}\\
\displaystyle{=z_2^{-1}\delta
\left(\frac{z_1-z_0}{z_2}\right)
Y(Y(u,z_0)v,z_2)}
\end{array}
\end{equation}
where 
$\delta(z)=\sum_{n\in {\Bbb Z}}z^n$ and  $(z_i-z_j)^n$ is 
expanded as a formal power series in $z_j.$ Throughout the paper,
$z_0,z_1,z_2,$ ect. are independent commuting formal variables.
\end{de}

As in [Z], we also define a  vertex operator superalgebra 
on torus associated to $V.$  Set $\tilde{\omega}=\omega-\frac{c}{24}$
and 
\begin{equation}\label{g2.8}
Y[v,z]=Y(v,e^z-1)e^{z{\rm wt}v}=\sum_{n\in\Z}v[n]z^{-n-1}
\end{equation}
for homogeneous $v\in V.$ Following the proof of Theorem 4.21 of \cite{Z} we have
\begin{th} $(V,Y[\ ],\bf{1},\tilde{\omega})$
is  a vertex operator superalgebra. 
\end{th}

Here we compute $v[m].$ Using a change of variable we calculate that 
\begin{eqnarray*}
& &v[m]=\Res_zY[v,z]z^m\\
& &\,\ \ \ \ \ \ =\Res_zY[v,\log(1+z)](\log(1+z))^m(1+z)^{-1}
\end{eqnarray*}
i.e.,
\begin{equation}\label{g2.12}
v[m]={\rm Res}_zY(v,z)(\log(1+z))^m(1+z)^{{\rm wt}v-1}.
\end{equation}

For integers $i,m,p$ with $i,m\geq 0,$ define numbers $c(p,i,m)$  by  
two equivalent ways:
\begin{equation}\label{g2.9}
{p-1+z\choose i}=\sum_{m=0}^ic(p,i,m)z^m.
\end {equation}
It is easy to see that
\begin {equation}\label{g2.10}
m!\sum_{i=m}^{\infty}c(p,i,m)z^i=(\log(1+z))^m(1+z)^{p-1}
\end{equation}
where, 
\begin{equation}\label{g2.11}
\log(1+z)=\sum_{n=1}^{\infty}\frac{(-1)^{n-1}}{n}z^n.
\end{equation}
If $m\geq 0$ then
\begin{equation}\label{g2.13}
v[m]=m!\sum_{i=m}^\infty c(\wt v,i,m)v(i).
\end{equation}
In particular,
\begin{equation}\label{g2.14}
v[0]=\sum_{i=0}^{\infty}{\wt v-1\choose i}v(i).
\end{equation}

Write
\begin{equation}\label{g2.15}
Y[\tilde \o,z]=\sum_{n\in\Z}L[n]z^{-n-2}.
\end{equation}
One can verify that
\begin{eqnarray}
& &L[-2]=\o[-1]-\frac{c}{24}\label{g2.16}\\
& & L[-1]=L(-1)+L(0)\label{g2.17}\\
& &\ \ \,L[0]=L(0)+\sum_{n=1}^{\infty}\frac{(-1)^{n-1}}{n(n+1)}L(n).\label{g2.18}
\end{eqnarray}
If $v\in V$ is homogeneous in  $(V,Y[\ ],{\bf 1},\tilde\o),$ 
we denote its conformal weight by $\wt [v].$ 

\begin{de}An {\rm automorphism} $g$ of vertex operator algebra $V$ is a
linear automorphism of $V$ preserving $\bf{1}$ and $\omega$ such that
the actions of $g$ and $Y(v,z)$ on $V$ are compatible in the sense
that $$gY(v,z)g^{-1}=Y(gv,z)$$ for $v\in V.$ 
\end{de}

Note that an automorphism preserves each homogeneous space $V_n$ and thus each $V_{\bar i}.$
There is a special automorphism $\sigma$ of $V$ with $\sigma|V_{\bar i}=(-1)^i$ associated to 
the superspace structure of $V.$ Then $\sigma$ is a central element 
in the full automorphism group $\Aut(V).$ As we will see later that $\sigma$ plays a fundamental
role in formulating the modular invariance of trace functions.   

\section{Twisted  Modules for VOSA}
\setcounter{equation}{0}

The twisted modules are the main objects in this section. In the
case of vertex operator algebras, the twisted modules are studied 
in \cite{D}, \cite{DLM2}, \cite{FFR}, \cite{FLM2}, \cite{FLM3},
\cite{L}. We will follow \cite{DLM2} and \cite{DZ} in this section.
In particular, we will introduce the notion of weak, admissible and
ordinary twisted modules. We also present 
the associative algebra $A_g(V)$ associated
to a vertex operator algebra $V$ and a finite order automorphism 
$g$ and related consequences.

Let  $(V,{\bf 1},\omega,Y)$ be a VOSA and $g$ an automorphism of $V$ of finite order $T.$
Let $T'$ be the order of $g\sigma.$ Then we have the following eigenspace
decompositions: 
 \begin{eqnarray}
& &V=\oplus_{r\in \Z/T'\Z}V^{r*} \label{g2.4}\\
& & V=\oplus_{r\in \Z/T\Z}V^{r} \label{g2.5} 				
\end{eqnarray}
where $V^{r*}=\{v\in V|g\sigma v=e^{-2\pi i r/T'}v\}$ and $V^{r}=\{v\in V|gv=e^{-2\pi ir/T}v\}.$
We now can define various notions of $g$-twisted modules (cf. \cite{DZ}).
\begin{de}A {\em weak} $g$-twisted $V$-module 
is a $\C$-linear space $M$ equipped with a linear map
$V\to (\End M)[[z^{1/T},z^{-1/T}]],$ $v\mapsto Y_M(v,z)=\sum_{n\in\Q}v(n) z^{-n-1}$ satisfying:

(1) For $v\in V, w\in M,$ $v(m) w=0$ if  $m >> 0,$

(2) $Y_M({\bf 1},z)=Id_{V},$

(3) For $v\in V^r,$ and $0\leq r\leq T-1,$ 

$$Y_M(v,z)=\sum_{n\in r/T+\Z}v(n)z^{-n-1}\ \ \ \ {\rm for}\ \ 
v\in V^r,$$

(4) The twisted Jacobi identity holds:
\begin{equation}\label{g3.6}
\begin{array}{c}
\displaystyle{z^{-1}_0\delta\left(\frac{z_1-z_2}{z_0}\right)
Y_M(u,z_1)Y_M(v,z_2)-\epsilon_{u,v}z^{-1}_0\delta\left(\frac{z_2-z_1}{-z_0}\right)
Y_M(v,z_2)Y_M(u,z_1)}\\
\displaystyle{=z_2^{-1}\left(\frac{z_1-z_0}{z_2}\right)^{-r/T}
\delta\left(\frac{z_1-z_0}{z_2}\right)
Y_M(Y(u,z_0)v,z_2)}
\end{array}
\end{equation}
if $u\in V^r.$
\end{de}

Note that (\ref{g3.6}) is equivalent to  the following:
\begin{eqnarray}\label{g3.7}
& &\Res_{z_1}Y_M(u,z_1)Y_M(v,z_2)z_1^{r/T}(z_1-z_2)^mz_2^n
\nonumber\\
& &\ \ \ -\epsilon_{u,v}\Res_{z_1}Y_M(v,z_2)Y_M(u,z_1)z_1^{r/T}(z_1-z_2)^m
z_2^n\nonumber\\
& &=\Res_{z_1-z_2}Y_M(Y(u,z_1-z_2)v,z_2)\iota_{z_2,z_1-z_2}z_1^{r/T}(z_1-z_2)^mz_2^n\label{ad2}
\end{eqnarray} 
for $m\in\Z, n\in\C,$ where $\iota_{z_2,z_1-z_2}z_1^{r/T}= \sum_{m\geq 0}{r/T\choose m}z_2^{r/T-m}
(z_1-z_2)^m.$ 

Set 
$$Y_M(\omega,z)=\sum_{n\in\Z}L(n)z^{-n-2}.$$
Then $Y_M(L(-1)z,v)=\frac{d}{dz}Y_M(v,z)$ for $v\in V$ and 
the operators $L(n)$ satisfy the Virasoro algebra relation with
central charge $c$ (cf. \cite{DLM1}, \cite{DZ}).

\begin{de}An {\em admissible} $g$-twisted $V$-module 
is a  weak $g$-twisted $V$-module $M$ which carries a 
${1\over T'}\Z$-grading 
\begin{equation}\label{m3.8}
M=\bigoplus_{0\leq n\in \frac{1}{T'}\Z}M(n)
\end{equation}
such that
\begin{eqnarray}\label{m3.9}
v(m)M(n)\subseteq M(n+\wt v-m-1)
\end{eqnarray}
for homogeneous $v\in V.$ We may and do assume that $M(0)\ne 0$ if $M\ne 0.$
\end{de}

Note that an admissible $V$-module $M$ is $\frac{1}{T'}\Z$-graded instead
of $\frac{1}{T}\Z$-graded as $V$ is $\frac{1}{2}\Z$-graded.

\begin{de}An (ordinary) $g$-twisted $V$-module $M$ is a $\C$-graded 
weak $g$-twisted $V$-module  
\begin{equation}\label{g3.10}
M=\coprod_{\lambda \in{\C}}M_{\lambda} 
\end{equation}
where $M_{\l}=\{w\in M|L(0)w=\l w\}$ such that
$\dim M_{\l}$ is finite and for fixed $\l,$ $M_{{n\over T}+\l}=0$
for all small enough integers $n.$ 
\end{de}

It $g=1$ we have the notion of weak, admissible ordinary $V$-modules.
It is not hard to prove that any ordinary $g$-twisted $V$-module is 
admissible.

If $M$ is a simple $g$-twisted $V$-module, then
\begin{equation}\label{g3.12}
M=\bigoplus_{n=0}^{\infty}M_{\la+n/T'}
\end{equation}
for some $\l\in\C$ such that $M_{\l}\ne 0$ (cf. \cite{Z}). 
We define $\l$ as the {\em conformal weight} of  $M.$

\begin{de} (1) A VOSA $V$ is called rational for an automorphism $g$ of finite order
if the category of admissible modules is semisimple.

(2) $V$ is called $g$-regular if any weak $g$-twisted $V$-module is a direct
sum of irreducible ordinary $g$-twisted $V$-modules.
\end{de}

It is clear that a $g$-regular vertex operator superalgebra is $g$-rational.

The following theorem is proved in \cite{DZ}.
\begin{th} Let $V$ be a VOSA and $g$ an automorphism of $V$ of finite
order. Assume that $V$ is $g$-rational. Then

(1) There are only finitely many irreducible admissible $g$-twisted $V$-modules
up to isomorphism.

(2) Each irreducible admissible $g$-twisted $V$-module is ordinary.
\end{th}

Now we review the $A_g(V)$ theory [DZ].  
Let $V^{r*}$ be as in (\ref {g2.4}).  For $u\in V^{r*}$ and $v\in V$ we define 
\begin{eqnarray}\label{g2.6}
u\circ_g v=\Res_{z}\frac{(1+z)^{{\wt u}-1+\delta_{r}+{r\over T'}}}{z^{1+\delta_{r}}}Y(u,z)v
\end{eqnarray}
\begin{equation}\label{g2.7}
u*_gv=\left\{
\begin{array}{ll}
\Res_z(Y(u,z)\frac{(1+z)^{{\wt}\,u}}{z}v)
 & {\rm if}\ r=0\\
0  & {\rm if}\ r>0.
\end{array}\right.
\end{equation} 
where $\delta_{r}=1$ if $r=0$ and 
$\delta_{r}=0$ if $r\ne 0$.
Extend $\circ_g$ and $*_g$ to bilinear products  on $V.$ We let $O_g(V)$ be the linear span of all $u\circ_g v.$
Then the quotient $A_g(V)=V/O_g(V)$ is an associative algebra with respect to 
$*_g.$

Let $M$ be a weak $\sigma g$-twisted $V$-module. Define $$\Omega (M)=\{w\in M | u(\wt u-1+n)w=0,u\in V,n>0\}.$$  For homogeneous $u\in V$ we set
\begin {equation}\label{zero mode}
o(u)=u(\wt u-1)
\end {equation}
We have the following Theorem (see \cite{DZ}):
\begin{th}\label{av} Let $V$ be a VOSA together with an automorphism $g$ of finite
order, and $M$ a weak $g$-twisted $V$-module. Then

(1) $\Omega(M)$ is an $A_g(V)$-module such that $v+O_g(V)$ acts as
$o(v).$

(2) If $M=\sum_{n\geq 0}M(n/T)$ is an admissible $g$-twisted $V$-module
with $M(0)\ne 0,$ then $M(0)\subset \Omega(M)$ is an $A_g(V)$-submodule.
Moreover, $M$ is irreducible if and only if $M(0)=\Omega(M)$ and
$M(0)$ is a simple $A_g(V)$-module.

(3) The map $M\to M(0)$ gives a 1-1 correspondence between the
irreducible admissible $g$-twisted $V$-modules and simple
$A_g(V)$-modules.

(4) If  $V$ is $g$-rational then $A_g(V)$ is a finite dimensional
semisimple associative algebra.
\end{th}
 
Another important concept is the $C_2$-cofinite condition (cf. \cite{Z}). Set
\begin{equation}\label{3.26}
C_2(V)=\{u(-2)v|u,\in V\}.
\end{equation}
$V$ is called $C_2$-cofinite if $C_2(V)$ has finite codimension.

It is proved in \cite{Li} and \cite{ABD} that a vertex operator algebra
$V$ is regular if and only if $V$ is $C_2$-cofinite and rational. Using the
the same proof we have
\begin{th} Assume that $V$ is $C_2$-cofinite. Then $V$ is  $g$-rational 
if and only if $V$ is $g$-regular.
\end{th}

Here is another consequence of the $C_2$-cofiniteness. 
\begin{prop}\label{p2.2} Suppose that $V$ satisfies condition $C_2,$
and let $g$ be an automorphism of $V$ of finite order. Then 

(1) The algebra $A_g(V)$ has finite dimension.

(2) If $A_g(V)\ne 0$ then $V$ has a simple $g$-twisted $V$-module.
\end{prop}
 
The proof of (1) is similar to that of Proposition 3.6 in \cite{DLM3} and
(2) follows from Theorem \ref{av} (3).

\section {{\bf 1}-point functions on the torus}
\setcounter{equation}{0}

We first introduce some notations which will be used in this section.
$V$ again is  a VOSA and $g,h\in Aut(V)$ such that $gh=hg.$
Let $o(g)=T, o(h)=T_1$  be finite.

Let $A$ be the subgroup of $Aut(V)$ generated by $g\sigma$ and
$h\sigma,$ and $N=lcm (T,T_1)$ be the exponent of $A.$ Let $\Gamma
(T,T_1)$ be the subgroup of matrices $\left(\begin{array}{cc} a & b\\
c & d
\end{array}\right)$ in $SL(2,\Z)$ satisfying $a\equiv d\equiv 1$ $(\mod\ N),$
$b\equiv 0 (\mod \ T),$ $c\equiv 0 (\mod \ T_1),$ and $M(T,T_1)$ be
the ring of holomorphic modular forms on $\Gamma (T,T_1)$ with natural
gradation $M(T,T_1)=\oplus_{k\geq 0}M_k(T,T_1),$ where $M_k(T,T_1)$ is
the space of forms of weight $k.$ From the Lemma 5.1 in [DLM3],
$M(T,T_1)$ is a Noetherian ring which contains each $E_{2k}(\t),$
$k\geq 2,$ and each $Q_k(\mu,\la,\t),$ $k\geq 0$ for $\mu,\la$ a
$T$-th., resp. $T_1$-th. root of unity.

Let $V(T,T_1)=M(T,T_1)\otimes_{\C}V.$ For $v\in V$ with $gv=\mu^{-1}v, hv=\la^{-1}v,$ we can define a $M(T,T_1)-$
submodule of $V(T,T_1),$ say $O(g,h)$ which consists of the following
elements:
\begin{eqnarray}
& &  v[0]w, w\in V, (\mu, \la)=(1,1)\label{m5.1}\\
& &  v[-2]w+\sum_{k=2}^{\infty}(2k-1)E_{2k}(\t)\otimes v[2k-2]w,  (\mu,\la)=(1,1)\label{m5.2}\\
& & v, (\tilde{v},\mu,\la)\ne (1,1,1)\label{m5.3}\\
& & \sum_{k=0}^{\infty} Q_k(\mu,\la,\t)\otimes v[k-1]w, (\mu,\la)\ne (1,1).
\label{m5.4}
\end{eqnarray}

We have the following Lemmas as in \cite{DLM3}.
\begin{lem} Suppose  $V$ satisfies $ C_2-$ cofinite condition. Then $ V(T,T_1)/O(g,h)$ is a finitely-generated $M(T,T_1)$-module.
\end{lem}
\begin{lem}\label{l5.3} If  $V$ is $C_2$-cofinite, then for  $v\in V,$  there is $m\in\N$ and elements $r_i(\tau)\in M(T,T_1),$ $0\leq i\leq m-1,$ 
such that
\begin{equation}\label{m5.6}
L[-2]^mv+\sum_{i=0}^{m-1}r_i(\tau)\otimes L[-2]^iv\in O(g,h).
\end{equation}
\end{lem}

\begin{de}\label{d5.3} The space of {\rm $(g,h)$ 1-point functions } 
${\cal C }(g,h)$ is a $\C$-linear space consisting of functions 

   $$S: V(T,T_1)\times {\frak h}\to \C $$  
satisfying  the following conditions:

(1) $S(v,\tau)$ is holomorphic in $\tau$ for $v\in V(T,T_1).$

(2) $S(v,\tau)$ is $\C$ linear in $v$ and  for $f\in M(T,T_1),$  $v\in V,$
 
 $$S(f\otimes v,\tau)=f(\tau)S(v,\tau)$$
 
(3) $S(v,\tau)=0$ if $v\in O(g,h).$

(4) If $v\in V$ with $\sigma v=gv=hv=v,$ then 

    \begin {equation}\label{m5.11a}
    S(L[-2]v,\tau)=\partial S(v,\tau)+\sum_{l=2}^{\infty}E_{2l}(\tau)S(L[2l-2]v,\t).
    \end{equation}
 Here  $ \partial S $ is the operator which is linear in $v$ and satisfies 
 
\begin{equation}\label{m5.11b}
\partial S(v,\t)=\partial_kS(v,\t)=\frac{1}{2\pi i}\frac{d}{d\t}S(v,\tau)
+kE_2(\tau)S(v,\t)
\end{equation}
for $v\in V_{[k]}.$
\end{de}

The following result -- modular transformation
is the main  property of $1$-point functions.

\begin{th}\label{t5.4}  For $S\in {\cal C}(g,h)$ and
$\gamma=\left(\begin{array}{cc}
a & b\\
c & d
\end{array}\right)\in \Gamma,$ we define 
\begin{equation}\label{m5.12}
S|\gamma(v,\t)=S|_k\gamma(v,\t)=(c\t+d)^{-k}S(v,\gamma\t)
\end{equation}
for $v\in V_{[k]},$ and extend linearly. Then $S|\gamma\in
{\cal C}((g,h)\gamma).$
\end{th}

This theorem is  the `super' analogue of Theorem 5.4 in [DLM3] and 
the proof is the same. 

Now we study the $q$-expansion of $S(v,\tau)$ for $S\in  {\cal C}(g,h)$
and $v\in V.$ As argued in \cite{DLM3}, the 
Frobenius-Fuchs theory tells us that $S(v,\tau)$ may be expressed
in the following form: for some $p\geq 0,$
\begin{equation}\label{6.12}
S(v,\tau)=\sum_{i=0}^p(\log q_{1\over T} )^iS_i(v,\tau)
\end{equation}
where
\begin{equation}\label{6.13}
S_i(v,\t)=\sum_{j=1}^{b(i)} q^{\lambda_{i,j}}S_{i,j}(v,\tau)
\end{equation}
\begin{equation}\label{6.14}
S_{i,j}(v,\tau)=\sum_{n=0}^{\infty}a_{i,j,n}(v)q^{n/T}
\end{equation}
are holomorphic on the upper half-plane, and
\begin{equation}\label{6.15}
\la_{i,j_1}\not\equiv \la_{i,j_2}\ (\mod \ \ \frac{1}{T}\Z)
\end{equation}
for $j_1\ne j_2$ where $q_{\frac{1}{T}}=e^{2\pi i\tau/T}$ and
$\la_{i,j}\in \C.$ 
 Moreover, $p$ is bounded independently of $v.$

It is important to point out that the 1-point functions on the torus
can be defined formally regarding $q$ as a formal variable. That is, 
we identify elements of $M(T,T_1)$ with
their Fourier expansions at $\infty,$ which lie in the ring
of formal power series $\C[[q_{1\over T} ]].$ Similarly, the functions $E_{2k}(\t),$
$k\geq 1,$ are considered to lie in $\C[[q]].$ 
The operator $\frac{1}{2\pi i}\frac{d}{d\tau}$ acts on $\C[[q_{1\over T} ]]$
via the identification $\frac{1}{2\pi i}\frac{d}{d \tau}=q\frac{d}{dq}.$

Then a  {\em formal} $(g,h)$ 1-point function is a map 
$$S: V(T,T_1)\to P$$
where $P$ is the space of formal power series of the form
\begin{equation}\label{7.1}
q^{\l}\sum_{n=0}^{\infty}a_nq^{n/T}
\end{equation}
for some $\l\in\C,$ and it satisfies the formal analogues
of (2)-(4) in the definition of $(g,h)$ 1-point functions. 
Then as in \cite{DLM3}, each $S(v,q)$ converges to a holomorphic
function. 

 we have 
\begin{th}\label{t7.1} Assume that $V$ is $C_2$-cofinite. Then any formal $(g,h)$ 1-point function $S$ defines an element of ${\cal C}(g,h),$ also
denoted by $S,$ via the identification
\begin{equation}\label{7.2}
S(v,\tau)=S(v,q), q=q_{\tau}=e^{2\pi i\tau}.
\end{equation}
\end{th}

\section{Trace functions}
\setcounter{equation}{0}
The goal in this section is to construct $(g,h)$ 1-point functions which
are essentially the graded $h\sigma$-trace function on $g\sigma$-twisted
$V$-modules.

Let $M=(M,Y_M)$ be a $g\sigma$-twisted $V$-module and $k\in \Aut (V).$
We can define  a $k(\sigma g) k^{-1}$-twisted $V$-module 
$(k\circ M,Y_{k\circ M})$ so that $k\circ M=M$ as vector spaces
and  
$$ Y_{k\circ M}(v,z)=Y_{M}(k^{-1}v,z).$$ 
$M$ is called $k$-stable if $k\circ M$ and $M$ are isomorphic.

Recall from Section 4 that if $M$ is a simple $g\sigma$-module then there 
exists a complex number $\lambda$ such that  
\begin{equation}\label{gg3.12}
M=\bigoplus^{\infty}_{n=0}M_{\lambda+\frac{n}{T}}
\end{equation}
(see (\ref{g3.12})) where $T$ is the order of $g.$ 
In particular, if $g=\sigma$ then 
\begin{equation}\label{gg3.12'}
M=\bigoplus^{\infty}_{n=0}M_{\lambda+\frac{n}{2}}.
\end{equation}

Let $h\in \Aut (V)$ as before and assume that $h\sigma\circ M$ and 
$M$ are isomorphic. Then $g,h$ commute and there is a linear
map $\phi(h\sigma): M\to M$ such that
\begin {equation}\label{6.1}
\phi(\sigma h)Y_M(v,z)\phi(\sigma h)^{-1}=Y_M((\sigma h)v ,z)    
\end{equation}  for all $v\in V.$ 

\begin{lem}\label{ll6.1} If $M$ is a simple $g\sigma$-twisted $V$-module then
$g\circ M$ and $M$ are isomorphic. In particular, any simple
$V$-module $M$ and $\sigma\circ M$ are isomorphic.
\end{lem}
\pf  We define a map $\phi(g)$ on $M$ such that 
$\phi(g)|_{M_{\l+n/T}}=e^{2\pi in/T}.$ Let $o(g\sigma)=T'$ and
$v\in V^{r*},$ $v\in V^{s}$  (see (\ref{g2.4})-(\ref{g2.5})). It is easy to verify that
$\wt v-\frac{r}{T'}$ and $\frac{s}{T}$ are congruent modulo $\Z.$ 
We immediately have 
$\phi(g)Y_M(v,z)\phi(g)^{-1}=Y_M(gv,z)$ for all $v\in V.$ 
That is, $M$ and $g\circ M$ are isomorphic. \qed

We remark that this result is different from a similar result in
the case of vertex operator algebra where a simple $g$-twisted module
$M$ and $g\circ M$ are always isomorphic. But this is not true in the
super case. 

Now we introduce the function $T$ which is linear in $v\in V,$ and define for homogeneous $v\in V$ as follows:
 \begin{equation}\label{6.4}
 T(v)=T_M(v,(g,h),q)=z^{\wt v} \tr_MY_M(v,z)\phi(\sigma h)q^{L(0)-\frac{c}{24}}
 \end{equation}
Here $c$ is the central charge of $V.$ Note that for $m\in
\frac{1}{T},v(m) $ maps $M_{\mu} $ to $M_{\mu+\wt v -m-1}.$ So we can
write $T(v)$ as follows:
\begin{equation}\label{6.5}
T(v)=q^{\lambda-\frac{c}{24}}\sum_{n=0}^{\infty}\tr_{M_{\lambda+\frac{n}{T}}}o(v)\phi(\sigma h)q^{\frac{n}{T}}=\tr_M o(v)\phi(\sigma h)q^{L(0)-\frac{c}{24}}
\end{equation} 

Here is the main result in this section.
\begin{th}\label{maint} Let $V$ be  $C_2$-cofinite and $g,h\in \Aut (V)$ have
finite orders. Let $M$  be a simple  $g\sigma$-twisted $V$-module 
such that $h\circ M, \sigma \circ M$ and $M$ are isomorphic. Then
$T\in {\cal C}(g,h).$
\end{th} 

By Theorem \ref{t7.1} it is enough to prove that  $T(v)$ defines a formal
$(g,h)$ 1-point function.  Clearly, $T(v)$ has the shape of (\ref{7.1}). 
So we need to prove that 
$T(v)$ satisfies (2)-(4) in definition \ref{d5.3}.
But (2) is clear,  we are going to prove (3) and (4). 

Fix $v\in V$ such that $gv=u^{-1}v,$ $hv=\lambda^{-1}v$ and  
$v_1\in V.$

\begin{lem}\label{t6.1}
 $T(v)=0$ for  $(\epsilon_{v},\mu,\la)\ne (1,1,1).$
\end{lem}

\pf Let $k\in \{g,h,\sigma\}.$ Since $k\circ M$ and $M$ are isomorphic
(see Lemma \ref{ll6.1}), then there is a linear isomorphism 
$\phi(k): M\to M$ such that $\phi(k)Y_M(v,z)\phi(k)^{-1}=Y_M(kv,z)$ 
for all $v\in V.$ 
Normalizing $\phi(k)$ we can assume
that $\phi(k)^t=1$ where $t$ is the order of $k.$ So we can decompose
$M$ into a direct sum of eigenspaces. Since $g,h,\sigma$ commute with
each other, $\phi(h\sigma)$ preserves each eigenspace. Since $v$ is 
an eigenvector for $k$ with eigenvalue different from 1, $o(v)$ moves
one eigenspace of $\phi(k)$ to a different eigenspace. As a result,
the trace of $o(v)\phi(h\sigma)$ on any homogeneous subspace
$M_{\l+n/T}$ is zero. The proof is complete.
\qed

\begin{lem}\label{t6.2}
$T(v[0]v_1)=0$  for  $(u,\lambda)=(1,1).$
\end{lem}
\pf Consider
    \begin {eqnarray*}
& & \ \    \tr_M[v(\wt v-1),Y_M(v_1,z)]\phi(\sigma h)q^{L(0)-c/24}\\
& &    =\!\tr_M v(\!\wt v\!-\!1\!)Y_M(v_1,z)\phi(\sigma h)q^{L(0)-c/24}\!-\!\epsilon_{v,v_1} \tr_M Y_M(v_1,z)v(\!\wt v\!-\!1\!
)\phi(\sigma h)q^{L(0)-c/24}\\
& &    =\!\tr_M v(\!\wt v\!-\!1\!)Y_M(v_1,z)\phi(\sigma h)q^{L(0)-c/24}\!-\!\epsilon_{v,v_1}\epsilon_{v}\tr_M Y_M(v_1,z)\phi(\sigma h)q^{L(0)-c/24}v(\!\wt v\!-\!1\!)\\ 
& &    =(1- \epsilon_{v,v_1}\epsilon_{v})\tr_M v(\wt v-1)Y_M(v_1,z)\phi(\sigma h)q^{L(0)-c/24}.
    \end{eqnarray*}
On the other hand
    \begin {eqnarray*}
  & & \ \  \tr_M[v(\wt v-1),Y_M(v_1,z)]\phi(\sigma h)q^{L(0)-c/24}\\
  & &  =\tr_M \sum^{\infty}_{i=0}{\wt v-1\choose i}z^{\wt v -1-i}Y_M(v(i)v_1,z)\phi(\sigma h)q^{L(0)-c/24}\\
  & &  =z^{-\wt v}T(v[0]v_1)      
    \end{eqnarray*} 
Thus if  $v\in V_{\bar{0}},$ or both $v$ and $v_1$ are in $V_{\bar{1}},$ 
then $\epsilon_{v,v_1}\epsilon_{v}=1$ and $T(v[0]v_1)=0.$ 
If $v\in V_{\bar{1}},$ and  
$v_1\in V_{\bar{0}},$ then $v[0]v_1\in V_{\bar{1}}.$ By Lemma \ref{t6.1}, 
we also have $T(v[0]v_1)=0.$\qed

\begin{lem}\label{t6.4}
 $$\sum_{k=0}^{\infty} Q_k(\mu,\la,\t)\otimes T(v[k-1]v_1)=0$$
 for $(\mu,\la)\ne (1,1).$
\end{lem}

In order to prove Lemma  \ref{t6.4}, we first define
$2$-point correlation functions. These are multi-linear functions $T(v_1,v_2)$ defined for $v_1,v_2 \in V$ homogeneous via
 \begin{eqnarray}\label{6.8}
& &T(v_1,v_2)=T((v_1,z_1),(v_2,z_2),(g,h),q)\nonumber\\
& &\ \ =z_1^{\wt v_1}z_2^{\wt v_2}{\tr}Y_M(v_1,z_1)Y_M(v_2,z_2)\phi(\sigma h)q^{L(0)-c/24}.
\end{eqnarray}       

We need the following sublemma. 
   
\begin{slem}\label{t6.7} Let $v,v_1\in V$ be homogeneous with $g v=\mu^{-1}v,$
$hv=\la^{-1}v$ and $(\mu,\la)\ne (1,1).$ Then
\begin{equation}\label{6.9a}
T(v,v_1)=\sum_{k=1}^{\infty}\bar P_{k}(\mu,\la,\frac{z_1}{z},q)T(v[k-1]v_1)
\end{equation}
\begin{equation}\label{6.9b}
T(v_1,v)=(-1)^{\tilde{v}}\la\sum_{k=1}^{\infty}\bar P_{k}(\mu,\la,\frac{z_1}{z}q,q)T(v[k-1]v_1)
\end {equation}
where $\bar P_k$ is as in (\ref{m4.29}).
\end{slem}  

\pf First we have the following identities:
\begin{eqnarray}\label{6.9c}
& &\ \ \  \ (1-\epsilon_{v,v_1}\epsilon_{v}\la q^k)\tr_M v(\wt v-1+k)Y_M(v_1,z_1)\phi(\sigma h)q^{L(0)}\nonumber\\
& &=\sum_{i=0}^{\infty}{\wt v-1+k\choose i}z_1^{\wt v-1+k-i}
\tr_M Y_M(v(i)v_1,z_1)\phi(\sigma h)q^{L(0)}
\end{eqnarray}
\begin{eqnarray}\label{6.9d}
& &\ \ \ \ (1-\epsilon_{v,v_1}\epsilon_{v}\la q^k)\tr_M Y_M(v_1,z_1)v(\wt v-1+k)\phi(\sigma h)q^{L(0)}\nonumber\\
& &=\epsilon_{v}\la q^k\sum_{i=0}^{\infty}{\wt v-1+k\choose i}z_1^{\wt v-1+k-i}
\tr_M Y_M(v(i)v_1,z_1)\phi(\sigma h)q^{L(0)}.
\end{eqnarray}
The proofs of these identities are similar to that given 
in Lemma 8.4 of \cite{DLM3} by noting the following
\begin{eqnarray}\label{6.10}
& &\tr_M v(\wt v-1+k)Y_M(v_1,z_1)\phi(\sigma h)q^{L(0)}=\tr_M [v(\wt v-1+k),Y_M(v_1,z_1)]\phi(\sigma h)q^{L(0)}\nonumber\\
& &\ \ \ \ \ +\epsilon_{v,v_1}\tr_M Y_M(v_1,z_1)v(\wt v-1+k)\phi(\sigma h)q^{L(0)},
\end{eqnarray}
$$v(\wt v-1+k)\phi(\sigma h)=\epsilon_{v}\la\phi(\sigma h)v(\wt v-1+k),$$
and 
\begin{equation}\label{6.11}
\tr Y_M(v_1,z_1)v(\wt v-1+k)\phi(\sigma h)q^{L(0)}
= \epsilon_{v}\la q^k\tr v(\wt v-1+k)Y_M(v_1,z_1)\phi(\sigma h)q^{L(0)}.
\end{equation}

Use identities (\ref{6.9c}), (\ref{6.9d}) and the same proof of Theorem
8.4 of \cite{DLM3} to obtain: 
\begin{eqnarray*}
&& T(v,v_1)=\sum_{k=1}^{\infty}\bar P_{k}(\mu,\epsilon_{v,v_1}\epsilon_{v}\la,\frac{z_1}{z},q)T(v[k-1]v_1)\\
&& T(v_1,v)=\epsilon_{v}\la^{-1}\sum_{k=1}^{\infty}\bar P_{k}(\mu,\epsilon_{v,v_1}\epsilon_{v}\la,\frac{z_1}{z}q,q)T(v[k-1]v_1).
\end{eqnarray*}
Again if $\epsilon_{v,v_1}\epsilon_{v}\ne 1$ then $T(v[k-1]v_1)=0$ by 
Lemma \ref{t6.1}. So in any case we have
\begin{eqnarray*}
&&T(v,v_1)=\sum_{k=1}^{\infty}\bar P_{k}(\mu,\la,\frac{z_1}{z},q)T(v[k-1]v_1)\\
&&T(v_1,v)=\epsilon_{v}\la^{-1}\sum_{k=1}^{\infty}\bar P_{k}(\mu,\la,\frac{z_1}{z}q,q)T(v[k-1]v_1),
\end{eqnarray*}
as required.\qed

We are now in a position to prove Lemma \ref{t6.4}. Recall the  following  
equality about Bernoulli polynomials from  \cite{DLM3}, Lemma 8.6.
\begin{equation}\label{eb}
 \sum_{k=0}^{\infty}\frac{1}{k!}B_k(1-\wt v+\frac{r}{T}+\frac{1}{2}\delta_{\tilde{v},1})v[k-1]=\sum_{i=0}^{\infty}
{\frac{r}{T}+\frac{1}{2}\delta_{\tilde{v},1}\choose i}v(i-1).
\end{equation}
  Let $m=\wt v-\frac{1}{2}\delta_{\tilde{v},1}\in \Z$  in Proposition 
\ref{p4.8}. Using (\ref{eb}) we have   
\begin{eqnarray*}
& &\ \ \sum_{k=0}^{\infty}Q_k(\mu,\la,q)T(v[k-1]v_1)\\
& &=\sum_{k=1}^{\infty}\Res_{z}\left((z-z_1)^{-1}z_1^{m -\frac{r}{T}-\frac{1}{2} \delta_{\tilde{v},1}}z^{-m+\frac{r}{T}+ \frac{1}{2}\delta_{\tilde{v},1}}\bar P_{k}(\mu,\la,\frac{z_1}{z},q)\right)T(v[k-1]v_1)\\
& &-\la \sum_{k=1}^{\infty} \Res_{z}\left((-z_1+z)^{-1}z_1^{m-\frac{r}{T}-\frac{1}{2} \delta_{\tilde{v},1}}z^{-m+\frac{r}{T}+ \frac{1}{2}\delta_{\tilde{v},1}}\bar P_{k}(\mu,\la,\frac{z_1q}{z},q)\right)T(v[k-1]v_1)\\
& & -\sum_{i=0}^{\infty}{\frac{r}{T}+\frac{1}{2} \delta_{\tilde{v},1}\choose i}T(v(i-1)v_1).
\end{eqnarray*}

On the other hand, by (\ref{g3.7}),
\begin{eqnarray*}
& &\sum_{i=0}^{\infty}{\frac{r}{T}+\frac{1}{2}\delta_{\tilde{v},1}\choose i}T(v(i-1)v_1)\\
& & =\sum_{i=0}^{\infty}{\frac{r}{T}+\frac{1}{2}\delta_{\tilde{v},1}\choose i}z_1^{\wt v+\wt v_1-i}\tr Y_M(v(i-1)v_1,z_1)\phi(\sigma h)q^{L(0)-c/24}\\
& & =\sum_{i=0}^{\infty}{\frac{r}{T}\!+\!\frac{1}{2}\delta_{\tilde{v},1}\choose i}z_1^{\wt v\!+\!\wt v_1\!-\!i}\Res_{z\!-\!z_1}(z\!-\!z_1)^{i\!-\!1}
\tr Y_M(Y(v,z\!-\!z_1)v_1,z_1)\phi(\sigma h)q^{L(0)\!-\!c/24}\\
 &&=\Res_{z\!-\!z_1}\iota_{z_1,z\!-\!z_1}\left(\frac{z}{z_1}\right)^{\frac{r}{T}\!+\!\frac{1}{2}\delta_{\tilde{v},1}}(\!z\!-\!z_1\!)^{-1}z_1^{\wt v\!+\!\wt v_1}
\tr Y_M(Y(v,z\!-\!z_1)v_1,z_1)\phi(\sigma h)q^{L(0)\!-\!c/24}\\
& &=\Res_z\iota_{z,z_1}(z-z_1)^{-1}({z_1\over z})^{\wt v-\frac{r}{T}-\frac{1}{2}\delta_{\tilde{v},1}}T(v,v_1)\\
& &\ \ \ -\epsilon_{v,v_1}\Res_z\iota_{z_1,z}(z-z_1)^{-1}({z_1\over z})^{\wt v-\frac{r}{T}-\frac{1}{2}\delta_{\tilde{v},1}}T(v_1,v)
\end{eqnarray*}
which by Sublemma \ref{t6.7} is equal to 
\begin{eqnarray*}
& &\sum_{k=1}^{\infty}\Res_z\iota_{z,z_1}(z-z_1)^{-1}({z_1\over z})^{\wt v-\frac{r}{T}-\frac{1}{2}\delta_{\tilde{v},1}}
\bar P_{k}(\mu,\la,\frac{z_1}{z},q)T(v[k-1]v_1)\\
& &-\epsilon_{v,v_1}\epsilon_{v}\la\sum_{k=1}^{\infty}\Res_z\iota_{z_1,z}(z-z_1)^{-1}({z_1\over z})^{\wt v-\frac{r}{T}-\frac{1}{2}\delta_{\tilde{v},1}}
\bar P_{k}(\mu,\la,\frac{z_1}{z}q,q)T(v[k-1]v_1).
\end{eqnarray*}
Again, if $v\in V_{\bar 0}$ or both $v,v_1\in V_{\bar 1}$ then
$\epsilon_{v,v_1}\epsilon_{v}=1$ and the result follows. If $v\in V_{\bar 1}$
and $v_1\in V_{\bar 0},$ the result follows from Lemma \ref{t6.1}. \qed

\begin{lem}\label{t6.3} If $(\mu,\la)=(1,1),$ then 
$$T(v[-2]v_1)+\sum_{k=2}^{\infty}(2k-1)E_{2k}(\t)T(v[2k-2]v_1)=0.$$
\end{lem}

\begin{slem}\label{t6.8} Let $v,v_1\in V$ be the same as in Lemma 
\ref{t6.3}.  Then
\begin{equation}\label{6.13.a}
T(v,v_1)=\tr_M o(v)o(v_1)\phi(\sigma h)q^{L(0)-\frac{c}{24}}+ \sum_{k=1}^{\infty}\bar P_{k}(1,1,\frac{z_1}{z},q)T(v[k-1]v_1)
\end{equation}
\begin{equation}\label{6.13.b}
T(v_1,v)=\tr_M o(v_1)o(v)\phi(\sigma h)q^{L(0)-\frac{c}{24}}+\epsilon_{v}\la\sum_{k=1}^{\infty}\bar P_{k}(1,1,\frac{z_1}{z}q,q)T(v[k-1]v_1)
\end {equation}
where $\bar P_k$ is defined  in (\ref{m4.29}).
\end{slem}

\pf We only prove (\ref{6.13.a}) and (\ref{6.13.b}) can be proved similarly.
By (\ref{6.9c}), (\ref{g2.9}) and (\ref{g2.13}) we have
\begin{eqnarray*}
& &T(v,v_1)=z^{\wt v}z_1^{\wt v_1}\tr Y_M(v,z)Y_M(v_1,z_1)\phi(\sigma h)q^{L(0)-\frac{c}{24}}\\
 & &\ \ \ =z^{\wt v}z_1^{\wt v_1}\tr \sum_{k\in \Z}z^{-k-\wt v}v(\wt v-1+k)Y_M(v_1,z_1)\phi(\sigma h)q^{L(0)-\frac{c}{24}}\\
&&\ \ \ =z_1^{\wt v_1}\tr o(v)Y_M(v_1,z_1)\phi(\sigma h)q^{L(0)-\frac{c}{24}}\\
&&\ \ \ +z_1^{\wt v_1}\sum_{k\in \Z- \{0\}}\frac{z^{-k}}{1\!-\!\epsilon_{v,v_1}\epsilon_{v}q^k}\sum_{i=0}^{\infty} {\wt v\!-\!1\!+\!k \choose i}z_1^{\wt v\!-\!1\!+\!k\!-\!i}\tr Y_M(v(i)v_1,z_1)\phi(\sigma h)q^{L(0)\!-\!\frac{c}{24}}\\
&&\ \ \ =z_1^{\wt v_1}\tr o(v)Y_M(v_1,z_1)\phi(\sigma h)q^{L(0)-\frac{c}{24}}\\
&&\ \ \ +\sum_{k\in \Z - \{0\}}|(\frac{z_1}{z})^k(1-\epsilon_{v,w}\epsilon_{v}q^k)^{-1}\sum_{i=0}^{\infty}\sum_{m=0}^{i}c(\wt v,i,m)k^mT(v(i)v_1)\\
& &\ \ \ = z_1^{\wt v_1}\tr o(v)Y_M(v_1,z_1)\phi(\sigma h)q^{L(0)-\frac{c}{24}}\\
 &&\ \ \ +\sum_{i=0}^{\infty}\sum_{m=0}^{i}\bar P_{m+1}(1,\epsilon_{v,v_1}\epsilon_{v},\frac{z_1}{z},q)m!c(\wt v,i,m)T(v(i)v_1)\\
& &\ \ \ =\tr o(v)o(v_1)\phi(\sigma h)q^{L(0)-c/24}+\sum_{m=0}^{\infty}\bar P_{m+1}(1,\epsilon_{v,v_1}\epsilon_{v},\frac{z_1}{z},q)T(v[m]v_1)\\
&&\ \ \ =\tr o(v)o(v_1)\phi(\sigma h)q^{L(0)-\frac{c}{24}}+\sum_{k=1}^{\infty}\bar P_{k}(1,\epsilon_{v,v_1}\epsilon_{v},\frac{z_1}{z},q)T(v[k-1]v_1) 
\end{eqnarray*}
Using Lemma \ref{t6.1} and discussing the values of $\epsilon_{v,v_1}\epsilon_{v}$  give the result. \qed

\begin{slem}\label{t6.8'} We have
 $$   T(v[-1]v_1)=\tr o(v)o(v_1)\phi(\sigma h)q^{L(0)-\frac{c}{24}}+\sum_{m=1}^{\infty}E_{2m}(\tau)T(v[2m-1]v_1).$$
   \end{slem}

\pf If  $\epsilon_{v,v_1}\epsilon_{v}\ne 1$ then both sides of the equation
are zero by Lemma \ref{t6.1}. So we now assume that  $\epsilon_{v,v_1}\epsilon_{v}=1.$ Write $v[-1]v_1=\sum_{i\geq -1}c_iv(i)v_1$ with $c_{-1}=1.$  Then 
\begin{eqnarray*}
&&T(v[-1]v_1)=\sum_{i\geq -1}c_iz_1^{\wt v+\wt v_1-i-1}\tr Y_M(v(i)v_1,z_1)\phi(\sigma h)q^{L(0)-\frac{c}{24}}\\
&&\ \ \ =\sum_{i\geq -1}c_iz_1^{\wt v+\wt v_1-i-1}\Res_{z_0}z_{0}^{i}\tr Y_M(Y(v,z_0)v_1,z_1)\phi(\sigma h)q^{L(0)-\frac{c}{24}}\\
&&\ \ \ =\sum_{i\geq -1}c_iz_1^{\wt v+\wt v_1-i-1}\Res_{z_{0}}z_{0}^{i}\Res_z(\frac{z_1+z_0}{z})^{\frac{1}{2}\delta_{\tilde{v},1}}\tr X\phi(\sigma h)q^{L(0)-\frac{c}{24}}
\end{eqnarray*}
where 
\begin{eqnarray*}
& & X=z_0^{-1}\delta(\frac{z-z_1}{z_0}) Y_M(v,z)Y_M(v_1,z_1)-\epsilon_{v,v_1}z_0^{-1}\delta(\frac{z_1-z}{-z_0}) Y_M(v_1,z_1)Y_M(v,z).
\end{eqnarray*}
Thus 
\begin{eqnarray*}
&&T(v[-1]v_1)=\sum_{i\geq -1}c_iz_1^{\wt v-i-1}\Res_z \sum_{j\geq 0}{\frac{1}{2}\delta_{\tilde{v},1} \choose j } z_1^{\frac{1}{2}\delta_{\tilde{v},1}-j}z^{-\frac{1}{2}\delta_{\tilde{v},1}-\wt  v}\cdot\\
&&\ \ \ \ \ \{(z-z_1)^{i+j}T(v,v_1)-\epsilon_{v,v_1}(-z_1+z)^{i+j}T(v_1,v)\}.
\end{eqnarray*}

Now replace $T(v,v_1)$ and $T(v_1,v)$ with (\ref{6.13.a}) and (\ref{6.13.b}) respectively and note that
$$\tr_M o(v_1)o(v)\phi(\sigma h)q^{L(0)-\frac{c}{24}}=
\epsilon_v\tr_M o(v)o(v_1)\phi(\sigma h)q^{L(0)-\frac{c}{24}}.$$
Then $T(v[-1]v_1)$ can be written as
\begin{equation}
 T(v[-1]v_1)=(a)+(b)
\end{equation}
where 
\begin{eqnarray*}
 & &(a)=\Res_z {z_{1}}^{\frac{1}{2}\delta_{\tilde{v},1}+\wt v}z^{ -\frac{1}{2}\delta_{\tilde{v},1}-\wt v}\{(z-z_1)^{-1}+(z_1-z)^{-1}\}\tr o(v)o(v_1)\phi(\sigma h)q^{L(0)-c/24}\\
     &&\ \ \ =\Res_z {z_{1}}^{\frac{1}{2}\delta_{\tilde{v},1}+\wt v}z^{ -\frac{1}{2}\delta_{\tilde{v},1}-\wt v}z_1^{-1}\delta\left(\frac{z}{z_1}\right)\tr o(v)o(v_1)\phi(\sigma h)q^{L(0)-\frac{c}{24}}\\
     &&\ \ \ =\tr o(v)o(v_1)\phi(\sigma h)q^{L(0)-\frac{c}{24}}
\end{eqnarray*}
(we have used the fact that $\frac{1}{2}\delta_{\tilde{v},1}+\wt v$ is always
an integer here)
and 
\begin{eqnarray*}
& &(b)=\sum_{k=1}^{\infty}\sum_{i\geq-1}c_i\Res_{z}\sum_{j\geq 0}{\frac{1}{2}\delta_{\tilde{v},1}\choose j}{z_1}^{\wt v+\frac{1}{2}\delta_{\tilde{v},1}-1-i-j}z^{-\frac{1}{2}\delta_{\tilde{v},1}-\wt v}\\
 && \ \ \ \ \{(z-z_1)^{i+j}P_{k}(1,1,\frac{z_1}{z},q)-(-z_1+z)^{i+j}P_{k}(1,1,\frac{z_1q}{z},q)\}T(v[k-1]v_1).
\end{eqnarray*}

From the proof of Proposition 4.3.3 of \cite{Z} we see that
$$(b)=\sum_{k=1}^{\infty}E_{2k}(q)T(v[2k-1]v_1).$$ 
This completes the proof.
\qed

Lemma \ref{t6.3} now follows from Sublemma \ref{t6.8'} (see the proof of Proposition 4.3.6 of \cite{Z}).

We also need the following Lemma whose  proof is contained in \cite{Z} by
 using Sublemma \ref{t6.8'}.
\begin{lem}\label{t6.5}
If $v\in V $ satisfies $ \sigma  v=gv=hv=v,$ then 
$$S(L[-2]v,\tau)=\partial S(v,\tau)+\sum_{l=2}^{\infty}E_{2l}(\tau)S(L[2l-2]v,\tau)$$
where
$$\partial S(v,\tau)=\partial_k S(v,\tau)=\frac{1}{2\pi i}\frac{d}{d\tau}S(v,\tau)+kE_2(\tau)S(L[2l-2]v,\tau) $$
\end{lem}

Theorem \ref{maint} follows from Lemmas \ref{t6.1}, \ref{t6.2},
\ref{t6.4}, \ref{t6.3} and \ref{t6.5}.

Here is another important theorem whose proof is similar to that of Theorem 
8.4 of \cite{DLM3}.       
        
 \begin{th}\label{t8.12} Let $M^1,M^2,...$ be inequivalent simple 
$\sigma g$-twisted
$V$-modules such that $h\sigma \circ M^i,$ $\sigma\circ M^i$ and
$M^i$ are isomorphic. Let $T_1,T_2,...$ be the
corresponding trace functions (\ref{6.4}). Then $T_1,T_2,...$ are linearly 
independent elements of ${\cal C}(g,h).$
\end{th}

\section{Existence of twisted modules}
\setcounter{equation}{0}

Although we have the definition of twisted module, but we do not know if there
is an irreducible $g$-twisted $V$-module for a finite order automorphism
$g$ of $V.$ The main result in this section gives a positive answer to
the question.

Let $g,h\in \Aut(V)$ commute and have finite orders. 
\begin{lem}\label{l9.3} Let $v\in V$ satisfy $g v=\mu^{-1}v, h v=\la^{-1}v.$
Then the following hold:

(1) The constant term of $\sum ^{\infty}_{k=0}Q_{k}(\mu,\la,q)v[k-1]w$ is equal
to $-v\circ_{\sigma g} w$ if $\mu\ne 1.$

(2) The constant term of $\sum_{k=0}^{\infty}Q_{k}(\mu,\la,q)(L[-1]v)[k-1]w$ is equal
to $-v\circ_{\sigma g} w$ if $\mu=1,$ $\la\ne 1.$

(3) The constant term of $v[-2]w+\sum_{k=2}^{\infty}(2k-1)E_{2k}(q)v[2k-2]w$ is 
$v\circ_{\sigma g} w$ if $\mu=\la=1.$
\end{lem}

\pf The proof is same as that of Lemma 9.2 in \cite{DLM3} by noting 
that $A_{\sigma g}(V)$ is defined by using  the decomposition of $V$ 
into a direct sum of eigenspaces  of $\sigma(\sigma g)=g.$
\qed 

\begin{th}\label{t9.1} Suppose that $V$ is a simple VOSA which satisfies 
condition $C_2,$ and that $g\in \Aut V$ has finite order. Then $V$ has 
at least one simple $g$-twisted module.
\end{th}

\pf Since $g$ is arbitrary, it is enough to prove that there exists a simple
$g\sigma$-twisted-module. By Proposition \ref{p2.2}, it suffices to
prove that $A_{\sigma g}(V)\ne 0.$ 

First we show that if ${\cal C}(g,h)\ne 0$ for some finite order 
$h\in \Aut(V)$ which commutes with $g$ then $A_{g\sigma}(V)$ is nonzero. 
As in \cite{DLM3} we take nonzero $S\in{\cal C}(g,h).$  Then
there exists positive integer $p$ such that (\ref{6.12})-(\ref{6.15}) hold for 
all $v\in V$ with $S_p\ne 0.$ We define a partial order on $\C$ as follows:
$\l\geq \mu $ if $\l-\mu\in\frac{n}{T}$ for some nonnegative integer
$\geq 0.$ We can arrange the notation
so that  $\l_{p,1}$ is  minimal among all $\l_{p,j}$ with respect to
the partial order and $a_{p,1,0}(v)\ne 0$ for some
$v\in V.$ Then
\begin{equation}\label{9.1}
S_{p,1}(v,\tau)=\a(v)+\sum_{n=1}^{\infty}a_{p,1,n}(v)q^{n/T}
\end{equation}
defines a function $\a: V\to\C$ which is not identically
zero. Because $S(v,\tau)$ is linear in $v,$ $\alpha$ is a linear
function on $V.$ Furthermore, $\alpha$ vanishes on $O_{\sigma g}(V)$
(see the proof of Lemma 9.3 of \cite{DLM3}). Thus $A_{g\sigma}(V)$ is nonzero.

Since $V$ is a simple $V$-module by the assumption,  
$T(v)=\tr_Vo(v)g\sigma q^{L(0)-c/24}$ defines a nonzero element in ${\cal C}(\sigma, g)$ by Theorems \ref{maint} and \ref{t8.12}. 
Since $\left(\begin{array}{cc} 0 & -1\\ 1 & 0\end{array}\right)$ induces 
a linear isomorphism between ${\cal C}(\sigma,g)$ and ${\cal C}(g,\sigma)$ by 
Theorem \ref{t5.4}. In particular, ${\cal C}(g,\sigma)$ is nonzero. The
proof is complete. \qed

\section{The main theorems}
\setcounter{equation}{0}

Recall Theorem \ref{t8.12}. In the case $V$ is $g\sigma$-rational, we
can strength the results obtained in \ref{t8.12}.
\begin{th}\label{t10.1}
Suppose that $V$ is $\sigma g$-rational and  $C_2$-cofinite.
Let $M^1,$ ..., $M^m$ be all of the inequivalent, simple
$\sigma g$-twisted $V$-modules such that $h\circ M^i,$ $\sigma \circ M^i$
and $M^i$ are isomorphic.  Let $T_1,...,T_m$ be the corresponding trace
functions defined by (\ref{6.4}). Then $T_1,...,T_m$ form a basis of ${\cal
C}(g,h).$ 
\end{th}

\begin{rem} In the proof of Theorem \ref{t10.1}, we only use
a consequence of $\sigma g$-rationality of $V$ --  $A_{\sigma g}(V)$
is finite dimensional semisimple algebra instead of $\sigma g$-rationality
of $V$ itself. So we can replace the $\sigma g$-rationality assumption
by the semi-simplicity of $A_{\sigma g}(V).$ 
\end{rem}

We begin with an arbitrary function $S\in{\cal C}( g, h).$ We have already
seen that $S$ can be represented as
\begin{equation}\label{10.1}
S(v,\tau)=\sum_{i=0}^p(log q_{1\over T} )^iS_i(v,\tau)
\end{equation}
for fixed $p$ and all $v\in V,$ with each $S_i$ satisfying (\ref{6.13})-(\ref{6.15}). Then Theorem \ref{t10.1} follows from the following propositions
\begin{prop}\label{p10.4} (1) Each $S_i$ is a linear combination of the functions $T_1,...,T_m.$

(2) $S_i=0$ if $i>0.$
\end{prop}

Using the proof of Proposition 10.5 in \cite{DLM3} we can show that
 Proposition (2) follows from (1). So we only need to prove (1).

We assume the setting in Section 7. In particular, 
$S_p\ne 0,$ each $S_{p,j}\ne 0$(\ref{6.13}), and  
$\a: V\to\C$ vanishes on $O_g(V)$ which 
induces a linear function
$$\alpha: A_{g\sigma}(V)\to \C.$$

\begin{lem}\label{l10.6}
Suppose that $u,v\in V$ and satisfy $hu=\rho u,hv=\la v,$ $\rho,\la \in \C.$ Then 
\begin{equation}\label{10.5}
\a(u*_{g\sigma} v)=\epsilon_{u,v}\rho\delta_{\rho\sigma,1}\a(v*_{g\sigma}u).
\end{equation}
\end{lem}

\pf Let $g u=\xi u$ and $g v=\nu v$ for scalars $\xi,\nu.$ If $\xi$ or $\nu$ is not equal to $1$ then $u$ (resp. $v$) lies in $O_{g\sigma}(V)$
by Lemma 3.1 of \cite{DZ}. This implies that
$u*_{g\sigma}v$ and $v*_{g\sigma}u$ are zero by Theorem 3.3 of \cite{DZ}.

We now assume $g u=u, g v=v.$ It is clear that
$u*_gv$ is an eigenvector for $h$ with eigenvalue $\rho\la.$ If $\rho\la\ne 1$ then $u*_gv$ and $v*_gu$ lie in $O(g,h)$ by (\ref{m5.3}).
Then $S(u*_{g\sigma}v)=S(v*_{g\sigma} u)=0$ by (C3). This leads to both sides 
of (\ref{10.5}) being $0.$ So finally we can assume that $\rho\la=1.$ 

Recall the decomposition (\ref{g2.4}). Then $u,v\in V^0.$ By  Lemma 3.3 of 
\cite{DZ}, 
$$ u*_{g\sigma}v-(-1)^{\tilde{u}\tilde{v}}v*_{g\sigma}u\equiv \Res_{z}Y(u,z)v(1+z)^{\wt u-1}\ (\mod  \ \ O(V^0)).$$
By (\ref{g2.14})
\begin{equation}\label{10.6}
 u*_{g\sigma}v-(-1)^{\tilde{u}\tilde{v}}v*_{g\sigma}u \equiv \sum_{i=0}^{\infty}{\wt u-1\choose i}u(i)v\equiv
u[0]v\ (\mod  \ \ O(V^0)).
\end{equation}
Since $O(V^{0})\subset O_{g\sigma}(V),$ 
$u[0]v\in O(g,h)$ if $\rho =1$  by (\ref{m5.1}). 
Thus $\a(u*_{g\sigma}v-\epsilon_{u,v}v*_{g\sigma}u)=0$ if $\rho=1.$

We now deal with the case $\rho \ne 1.$ From the proof of Lemma 9.2
of \cite{DLM3}, we see that the constant term  of 
$$\sum_{k=0}^{\infty}Q_k(1,\rho^{-1},q)u[k-1]v$$
is 
$$-u*_{g\sigma}v+\frac{1}{1-\rho^{-1}}u[0]v.$$
The vanishing condition of $S$  on $\sum_{k=0}^{\infty}Q_k(1,\rho^{-1},q)u[k-1]v\in O(g,h)$ yields  
$$\a(u*_g v)=\frac{1}{1-\rho^{-1}}\a(u[0]v).$$
Using (\ref{10.6}) gives the desired result. \qed

\begin{lem}\label{*} Suppose that $u,v\in V$ with  $h\sigma u=\rho' u,h\sigma v=\la' v,$ $\rho',\la' \in \C.$ Then 
\begin{equation}\label{**}
\a(u*_{g\sigma} v)=\rho'\delta_{\rho'\la',1}\a(v*_{g\sigma}u).
\end{equation}
\end{lem}

\pf We still assume that  $hu=\rho u,hv=\la v.$ Then $\rho'=\epsilon_u\rho$ and
$\la'=\epsilon_v\la.$ By Lemma \ref{l10.6} we have
$$\a(u*_{g\sigma} v)=\epsilon_{u,v} \epsilon_{u}\rho'\delta_{\epsilon_{u}\epsilon_{v}\rho'\la',1}\a(v*_{g\sigma}u).$$

If $\epsilon_{u}\epsilon_{v}=-1,$ then  
$\sigma(u*_{g\sigma}v)=-u*_{g\sigma}v.$ By (\ref{m5.3}),both $u*_{g\sigma}v$ and $v*_{g\sigma} u$  lie in 
$O(g,h).$ This implies that $S(u*_{g\sigma}v)=S(v*_{g\sigma} u)=0.$ As 
a result, both sides of (\ref{**}) are equal to $0.$
If $\epsilon_{u}\epsilon_{v}=1,$ then $\epsilon_{u,v}\epsilon_{u}=1$ and
(\ref{**}) still holds. \qed

The rest of the proof of Proposition \ref{p10.4} (1) follows from 
the argument given in Section 10 of \cite{DLM3}.

As in \cite{DLM3} we have important corollaries.
\begin{th}\label{t10.2}  Suppose that $V$ is rational and $C_2$-cofinite. 

(1) If the group $\<g,h\>$ generated by $g$ and $h$ is cyclic
with generator $k.$ Then the dimension of ${\cal C}( g, h)$ is equal
to the number of inequivalent, $k,\sigma$-stable, $\sigma$-twisted
simple $V$-modules. 

(2) The number of inequivalent, $\sigma$-stable simple $g$-twisted 
$V$-modules is at most equal to the number of $g,\sigma$-stable, simple $V$-modules,
with equality if $V$ is $g$-rational.
\end{th}

The proof is the same as that of Theorem 10.2 in \cite{DLM3} by using Lemma
\ref{ll6.1} and  Theorem \ref{t10.1}.

As in the theory of vertex operator algebra, a 
simple vertex operator superalgebra $V$ is called {\em holomorphic} in 
case $V$ 
is rational and if $V$ is the unique simple $V$-module.

\begin{th}\label{t10.3} Suppose that $V$ is holomorphic and is  
$C_2$-cofinite.  For each automorphism $g$ of $V$ of finite order, there is
a unique $\sigma$-stable simple $g$-twisted $V$-module $V(g).$ Moreover 
if $\<g,h\>$
is cyclic then ${\cal C }( g, h)$ is spanned by $T_{V(\sigma g)}(v,g,h,q).$
\end{th}

It is proved in \cite{Z} that if a vertex operator algebra
$V$ is rational and $C_2$-cofinite, then the space
spanned by $T_M(v,q)=\tr_Mo(v)q^{L(0)-c/24}$ with $M$ running through
the inequivalent simple $V$-modules admits a representation
of the modular group $\Gamma.$ But this is not true anymore 
in the present situation. In fact, if $V$ is holomorphic such 
that $V_{\bar 1}\ne 0,$ then the character $\tr_Vq^{L(0)-c/24}$ of $V$
is not an eigenvector for the matrix
$T=\left(\begin{array}{cc}1 & 1\\ 0 &1\end{array}\right).$
Here is the corresponding result for the vertex operator
superalgebra, which is an immediate consequence of Theorems
\ref{maint} and \ref{t10.1}.
\begin{th} Let $V$ be $\sigma$-rational and $C_2$-cofinite
VOSA. Let $M^1,$ ..., $M^m$ be all of the inequivalent, simple
$\sigma $-twisted $V$-modules such that $\sigma \circ M^i$ 
and $M^i$ are isomorphic. Then the space
spanned by 
$$T_i(v,\tau)=T_i(v,(1,1),\tau)=\tr_{M^i}o(v)\phi(\sigma)q^{L(0)-c/24}$$admits a representation of the modular group. That is, for 
any $\gamma=\left(\begin{array}{cc}a & b\\ c &d\end{array}\right)
\in \Gamma$ there exists a  $m\times m$ invertible 
complex matrices $(\gamma_{ij})$ such that
$$T_i(v,\frac{a\tau+b}{c\tau+d})=(c\tau+d)^n\sum_{j=1}^m\gamma_{ij}T_j(v,\tau)$$for all $v\in V_{[n]}.$ Moreover, the matrix $(\gamma_{ij})$ is independent
of $v.$
\end{th}

Finally we discuss the rationality of the central charge $c$ and 
the conformal weights of a $g$-twisted module under certain conditions.
Let $M$ be a simple $g$-twisted $V$-module. Recall from (\ref{g3.12})
that $M$ is direct sum of $L(0)$-eigenspaces 
$$M=\bigoplus_{n=0}^{\infty}M_{\la+n/T'}$$
for some $\l\in\C$ such that $M_{\l}\ne 0.$ The $\l$ is 
called the conformal weight of  $M.$ 

\begin{th}\label{t11.2} Assume that  $V$ is $C_2$-cofinite and 
$g\in\Aut V$ has finite order. Suppose that $V$ is
$\sigma^i g^j$-rational for all integers $i,j.$ Then each  $\sigma$-stable
simple $\sigma^i g^j$-twisted $V$-module has rational conformal weight, and the central 
charge $c$ of $V$ is rational. 
\end{th}

The proof of this theorem is the same as that of Theorem 11.2 in \cite{DLM3}
(also see \cite{AM}) using Lemma \ref{ll6.1}.

\section{Modular invariance over $\Gamma_{\theta}$}

In order to get the modular invariance over the full modular group $\Gamma$ 
we only considered the $\sigma$-stable twisted modules in previous sections.
In this section we remove the $\sigma$-stable assumption. But we 
have a modular invariance over the $\Gamma_{\theta}$ instead of 
the full modular group $\Gamma$ as in \cite{H}. Since the proofs of most
 results in this section are very similar to those given in the previous 
sections and \cite{DLM3} we present the results without proofs.

Recall that $G$ is a finite automorphism group of $V.$ Let $g\in G$ and
$M$ a simple $g$-twisted $V$-module. Also recall $M\circ h$ for
$h\in G.$ Let $G_M=\{h\in G|M\circ h\cong M\}$ be the stablizer of $M.$
Then $M$ is a projective module for $G_M,$ denoting by $\phi$ the 
projective action. For $h\in G_M$ and $v\in V$ set
$$F_M(v, (g, h),q)  = \tr_M o(v)\phi(h)q^{L(0)-c/24}.$$

\begin{th} Let $V$ be a rational and $C_2$-cofinite vertex operator superalgebra.

(1)  If $g$ is an automorphism of $V$ of finite order then the number of
inequivalent, $g$-stable, simple  $g$-twisted $V$-modules is at most equal to 
the number of $g$-stable simple untwisted $V$-modules.
 If $V$ is $g$-rational, the number of
inequivalent, $g$-stable simple $g$-twisted $V$-modules is precisely the
number of $g$-stable simple untwisted $V$-modules.

(2) Let $g\in\Aut V$ have finite order. Suppose that $V$ is
$g^i$-rational for all integers $i.$ Then each $g$-stable simple $g^i$-twisted
$V$-module has rational conformal weight, and the central charge $c$ of $V$
is rational. 

(3) Let $G$ be a finite automorphism group of $V.$ Let $g\in G$ and $M$ a simple $g$-twisted $V$-module
and $h\in G_M.$ Then 
the trace function  $F_M(v,(g,h),q)$ converges
to a holomorphic function in the upper half plane ${\frak h}.$ 

(4) Assume that $V$ is $x$-rational for each $x\in G.$ 
Let $v\in V$ satisfy $\wt[v] = k.$ Then the space of (holomorphic) functions 
in ${\frak h}$ spanned
by the trace functions $F_M(v,(g,h), \tau)$ for all choices of $g, h$ in $G$
 and 
$h$-stable $M$  is  a (finite-dimensional)  $\Gamma_{\theta}$-module. That is,
if $\gamma=\left(\begin{array}{cc}
a & b\\ c& d\end{array}\right)\in \Gamma_{\theta}$ then we have an equality
$$F_M(v, (g, h),\frac{a\tau+b}{c\tau+d})=(c\tau+d)^k\sum_{W}\gamma_{M,W}F_W(v,(g^ah^c,g^bh^d),\tau),$$
 \end{th}

The statements (3) and (4) are exactly the same as in \cite{DLM3} for
vertex operator algebras except that we use the group $\Gamma_{\theta}$
here instead of $\Gamma$ in \cite{DLM3}. In the 
case $G=1,$ Parts (3) and (4) have been obtained previously in \cite{H}.

\section {Examples}

 In this section we will use the examples in [DZ] (cf. Theorem
   7.1, propositions 7.2 and 7.4) to show the modular
   invariance. So we are working in the setting of
Section 7 of \cite{DZ}. In particular, $H$ is a vector space 
of dimension $l$ which we assume to be even, $V(H,{\Z}+\frac{1}{2})$ is a 
VOSA with central charge $c=\frac{l}{2}$ and $V(H,\Z)$ is the unique 
irreducible $\sigma$-twisted module. We are going to study the 
trace functions $T({\1}, (\sigma^i,\sigma^j),q)$ which is the graded 
 $\sigma^{j+1}-$trace on the irreducible $\sigma^{i+1}$-twisted module
for $V(H,{\Z}+\frac{1}{2})$ for $i,j=0,1.$ 

We first recall the important modular form  
$$\eta(\tau)=q^{\frac{1}{24}}\prod_{n\geq 1}(1-q^n).$$
The trace functions in our examples can be expressed as rational functions 
in $\eta(\tau).$ Let 
$S=
\left(\begin{array}{cc}
0 & 1\\
-1 & 0
\end{array}
\right)$  and $T=\left(\begin{array}{cc}
1 & 1\\
0 & 1
\end{array}
\right)$ be the standard generators of $SL(2,\Z).$ Then 
$$\eta(\frac{-1}{\tau})=(-i\tau)^{\frac{1}{2}}\eta(\tau)$$  
and 
$$\eta(\tau+1)=e^{\frac{\pi i}{12}}\eta(\tau).$$

Note that the conformal weight of $V(H,\Z)$ is $\lambda=\frac{l}{16}.$ It is
easy to compute that 
\begin{eqnarray*}
& &T({\1}, (1,1),\tau)=\tr_{V(H,\Z)}\sigma q^{L(0)-\frac{l}{48}}=q^{\frac{l}{24}}\prod_{m=1}^{l/2}(1+(-1))\prod_{n>0}(1-q^n)^l=0\\
& & T(\1,(1,\sigma),\tau)=\tr_{V(H,\Z)}q^{L(0)-\frac{l}{48}}
=2^{l/2}q^{\frac{l}{24}}\prod_{n=1}^{\infty}(1+q^n)^l=2^{l/2}[\frac{\eta(2\tau)}{\eta(\tau)}]^l\\
& &T(\1,(\sigma, 1),\tau)= \tr_{V(H,\frac{1}{2}+\Z)}\sigma q^{L(0)-\frac{l}{48}}=q^{-\frac{l}{48}}\prod_{n=0}^{\infty}(1-q^{n+\frac{1}{2}})^l=[\frac{\eta(\tau/2)}{\eta(\tau)}]^l\\
& & T(\1,(\sigma, \sigma),\tau)= \tr_{V(H,\frac{1}{2}+\Z)}q^{L(0)-\frac{l}{48}}
=q^{-\frac{l}{48}}\prod_{n=0}^{\infty}(1+q^{n+\frac{1}{2}})^l
=[\frac{\eta(\tau)^2}{\eta(\tau/2)\eta(2\tau)}]^l.
 \end{eqnarray*} 
The modular transformation law then gives the modular transformation
property of functions $T(\1, (\sigma^i,\sigma^j),\tau)$ for
$i,j=0,1:$ 
$$T(\1, (\sigma^i,\sigma^j),\frac{a\tau+b}{c\tau+d})=\gamma_{i,j}T(\1, (\sigma^{ai+cj},\sigma^{bi+dj}),\tau)$$
for some constant $\gamma_{i,j}\in\C$ where $\gamma=\left(\begin{array}{cc}
a & b\\
c & d
\end{array}
\right)\in SL(2,\Z).$  
The result matches what Theorems \ref{t5.4} and \ref{t10.1} assert.
     
Here is another example. We still take $V$ to be $V(H,\frac{1}{2}+\Z)$ with
$l=\dim H$ being even. Let $g$ be an isometry of $H$ such that
$g(h_1)=-h_2, g(h_2)=-h_1,$ $g(h_1^*)=-h_2^*, g(h_2^*)=-h_1^*.$ And $ g(h_i)=-h_i,g(h_i^*)=-h_i^*$ for $i=3,\cdots l/2$ where
$\{h_1,...,h_{l/2}, h_1^*,...,h_{l/2}^*\}$ is a basis of $H$ such that
$(h_i,h_j)=(h_i^*,h_j^*)=0$ and $(h_i,h_j^*)=\delta_{i,j}.$ 
Then $g\sigma$ interchanges $h_1$ and $h_2,$ $h_1^*$ and $h^*_2,$ and
leaves other $h_i,$ $h_i^*$ invariant. Then 
we have 
$$H=H^{0*}\bigoplus H^{1*}, \ \ \  H=H^0\bigoplus H^1$$
where $g\sigma|_{H^{i*}}=(-1)^i$ and $g|_{H^i}=(-1)^i.$ It is easy
to see that 
$$H^{0*}=\C(h_1+h_2)\bigoplus\C(h_1^*+h_2^*)\bigoplus\bigoplus_{j=3}^{l-2}(\C h_j+\C h_j^*),$$
$$H^{1*}=\C(h_1-h_2)\bigoplus\C(h_1^*-h_2^*),$$ 
and $H^0=H^{1*},H^1=H^{0*}.$ Note that $g$ is lifted to an automorphism 
of the vertex operator superalgebra $V(H,\Z+\frac{1}{2})$ in an obvious
way.

By proposition 7.4 in \cite{DZ},
$$M=\wedge[(h_1-h_2)(-n), (h_1^*-h_2^*)(-m),h(-m-\frac{1}{2}) |n,m\in \Z,n>0,m\geq0, h\in H^1]$$
 is the unique irreducible $g\sigma$-twisted $V(H,\Z+\frac{1}{2})$-module 
whose conformal weight  is $\frac{1}{8}.$
Thus
\begin{eqnarray*}
&&T(\1,(g,\sigma),\tau)=\tr_Mq^{L(0)-\frac{l}{48}}\\
&&=\sum_{n\in \frac{1}{2}\Z}(\dim M_{n+\frac{1}{8}})q^{n-\frac{l}{48}}\\\nonumber
&&=2[\prod_{0<n\in \Z}(1+q^n)]^2 [\prod_{0\leq n\in \Z}(1+q^{n+\frac{1}{2}})]^{l-2}q^{\frac{1}{8}-\frac{l}{48}}\\\nonumber
&&=2[q^{-\frac{1}{24}}\frac{\eta(2\tau)}{\eta(\tau)}]^2[q^{\frac{1}{48}}\frac{\eta^2(\tau)}{\eta(\frac{\tau}{2})\eta(2\tau)}]^{l-2}q^{\frac{1}{8}-\frac{l}{48}}\\\nonumber
&&=2\frac{\eta(\tau)^{2l-6}}{\eta(2\tau)^{l-4}\eta(\frac{\tau}{2})^{l-2}},
\end{eqnarray*}
\begin{eqnarray*}
&&T(\1,(\sigma,g),\tau)=\tr_{V(H,\Z+\frac{1}{2})}g\sigma q^{L(0)-\frac{l}{48}}\\
&&=\sum_{0\leq n\in \frac{1}{2}{\Z}}(\tr_{ V(H,\Z+\frac{1}{2})_n}g\sigma)q^{n-\frac{l}{48}}\\\nonumber
&&=[\prod_{0\leq n\in \Z}(1-q^{n+\frac{1}{2}})]^2 [\prod_{0\leq n \in \Z}(1+q^{n+\frac{1}{2}})]^{l-2}q^{-\frac{l}{48}}\\\nonumber
&&=[q^{\frac{1}{48}}\frac{\eta(\frac{\tau}{2})}{\eta(\tau)}]^2[q^{\frac{l}{48}}\frac{\eta^2(\tau)}{\eta(\frac{\tau}{2})\eta(2\tau)}]^{l-2}q^{-l/48}\\\nonumber
&&=\frac{\eta(\tau)^{2l-6}}{\eta(\frac{\tau}{2})^{l-4}\eta(2\tau)^{l-2}}.
\end{eqnarray*}
and 
\begin{eqnarray*}
&&T(\1,(g,g\sigma),\tau)=\tr_Mgq^{L(0)-\frac{l}{48}}\\
& &\ \ \ =\sum_{0\leq n\in \frac{1}{2}{\Z}}(tr_{ M_{n+\frac{1}{8}}}g)q^{L(0)-\frac{l}{48}}\\\nonumber
&&\ \ \ =2q^{\frac{1}{8}-\frac{l}{48}}[\prod_{0<n\in \Z}(1+q^n)]^2[\prod_{0\leq n\in \Z}(1-q^{n+\frac{1}{2}})]^{l-2}\\\nonumber
&&\ \ \ =2q^{\frac{1}{8}-\frac{l}{48}}[q^{\frac{-1}{24}}\frac{\eta(2\tau)}{\eta(\tau)}]^2[q^{\frac{1}{48}}\frac{\eta(\frac{\tau}{2})}{\eta(\tau)}]^{l-2}\\\nonumber
&&=\ \ \ 2\frac{\eta(2\tau)^2\eta(\frac{\tau}{2})^{l-2}}{\eta(\tau)^l}
\end{eqnarray*}

Then 
\begin {eqnarray*}
&&T(\1,(g,\sigma),S\tau)=2\frac{\eta(-\frac{1}{\tau})^{2l-6}}{\eta(-\frac{2}{\tau})^{l-4}\eta(\frac{-1}{2\tau})^{l-2}}\\\nonumber
&&\ \ \ =2\frac{(-i\tau)^{\frac{2l-6}{2}}\eta(\tau)^{2l-6}}{(-i\tau/2)^{\frac{l-4}{2}}\eta(\frac{\tau}{2})^{l-4}(-i2\tau)^{\frac{l-2}{2}}\eta(2\tau)^{l-2}}\\\nonumber
&&\ \ \ =2\mu \frac{\eta(\tau)^{2l-6}}{\eta(\frac{\tau}{2})^{l-4}\eta(2\tau)^{l-2}}
\end{eqnarray*}    
where $\mu$ is a root of unity. This implies that 
$$T(\1,(g,\sigma),S\tau)=2\mu T(\1, (g,\sigma)S,\tau)=2\mu T(\1, (\sigma,g),\tau).$$

Using the fact that $$\eta(\frac{\tau+1}{2})=\eta(\tau)\frac{\eta(\tau)^2}{\eta(\frac{\tau}{2})\eta(2\tau)}$$
we have 
\begin{eqnarray*}
&&T(\1,(g,\sigma),T\tau)=T(\1,(g,\sigma),\tau+1)\\\nonumber
&&\ \ \ =2\frac{\eta(\tau+1)^{2l-6}}{\eta(2(\tau+1))^{l-4}\eta(\frac{\tau+1}{2})^{l-2}}\\\nonumber
&&\ \ \ =2\nu \frac{\eta(\tau)^{2l-6}}{\eta(2\tau)^{l-4}}\frac{\eta(\frac{\tau}{2})^{l-2}\eta(2\tau)^{l-2}}{\eta(\tau)^{2l-4}}\\\nonumber
&&\ \ \ =2\nu \frac{\eta(2\tau)^2\eta(\frac{\tau}{2})^{l-2}}{\eta(\tau)^l}
\end{eqnarray*}
for some root of unity $\nu.$ Since $(g,\sigma)T=(g,g\sigma),$ we see
that
$$T(\1,(g,\sigma),T\tau)=\nu T(\1,(g,g\sigma),\tau).$$
This again verifies Theorems \ref{t5.4} and \ref{t10.1} in this special
case.  The reader could determine $T(v,(h,k),\tau)$ for other  
vector $v\in V$ in these two cases, but the computation will be
more complicated and difficult.

\end{document}